\documentclass[11pt,oneside]{amsart}
\usepackage{amscd,amssymb}
\usepackage{graphicx}

\theoremstyle{plain}
\newtheorem{Thm}[equation]{Theorem}

\newtheorem{Rmk}[equation]{Remark}
\newtheorem{Cor}[equation]{Corollary}
\newtheorem{Prop}[equation]{Proposition}

\newtheorem{Lem}[equation]{Lemma}
\newtheorem{Def}[equation]{Definition}
\numberwithin{equation}{section}

\newcommand{\e}{\epsilon}
\newcommand{\z}{\mathbb{Z}}

\renewcommand{\c}{\mathbb{C}}
\newcommand{\br}{\mathbb{R}}

\newcommand{\A}{\mathcal{A}}
\newcommand{\ba}{\backslash}

\newcommand{\G}{\Gamma}

\renewcommand{\P}{\mathcal P}

\newcommand{\la}{\langle}
\newcommand{\ra}{\rangle}

\newcommand{\bp}{\begin{pmatrix}}
\newcommand{\ep}{\end{pmatrix}}

\newcommand{\PSL}{\op{PSL}}
\newcommand{\BMS}{\op{BMS}}
\newcommand{\BR}{\op{BR}}
\newcommand{\bi}{\begin{itemize}}
\newcommand{\be}{\begin{enumerate}}

\newcommand{\ee}{\end{enumerate}}

\newcommand{\op}{\operatorname}

\newcommand{\vs}{\vskip 5pt}

\newcommand{\PS}{\operatorname{PS}}
\newcommand{\sk}{\operatorname{sk}}
\newcommand{\Leb}{\operatorname{Leb}}
\newcommand{\bH}{\mathbb H}
\newcommand{\T}{\op{T}}
\renewcommand{\PS}{\operatorname{PS}}
\renewcommand{\BR}{\operatorname{BR}}
\renewcommand{\Leb}{\operatorname{Leb}}
\renewcommand{\BMS}{\operatorname{BMS}}

\begin{document}

\title[Asymptotic distribution of circles in orbits of Kleinian groups]{The asymptotic distribution of circles in the orbits of Kleinian groups}

\author{Hee Oh and Nimish Shah}

\address{Mathematics department, Brown university, Providence, RI
and Korea Institute for Advanced Study, Seoul, Korea}
\email{heeoh@math.brown.edu}

\address{Department of Mathematics, The Ohio State University, Columbus, OH}
\email{shah@math.ohio-state.edu}

\thanks{Oh is partially supported by NSF
   grant 0629322}
\thanks{Shah is partially supported by NSF grant 1001654}
\begin{abstract} Let $\P$ be a locally finite circle packing in the plane $\c$ invariant under
a non-elementary Kleinian group $\G$ and with finitely many $\G$-orbits.
When $\G$ is geometrically finite,
we construct an explicit Borel measure on $\c$ which
describes the asymptotic distribution of small circles in $\P$, assuming that
either the critical exponent of $\G$ is strictly bigger than $1$ or
$\P$ does not contain an infinite bouquet of tangent circles glued at a parabolic fixed point of $\G$.
Our construction also works for $\P$ invariant under a
geometrically {\it infinite} group $\G$, provided
 $\G$ admits a finite Bowen-Margulis-Sullivan measure and the $\G$-skinning size of $\P$ is finite.
Some concrete circle packings to which our result applies include Apollonian circle packings, Sierpinski curves,
 Schottky dances, etc.
\end{abstract}

\maketitle

\section{Introduction}
A circle packing in the plane $\c$ is simply a union of circles (here a line
is regarded as a circle of infinite radius). As we allow circles to intersect with each other,
our definition of a circle packing is more general than the conventional definition of a 
circle packing. 
 
 For a given circle packing $\P$ in the plane, we are interested in counting and distribution
 of small circles in $\P$. A natural size of a circle is measured by its radius. We will use
the curvature of a circle, that is, the reciprocal of its radius, instead.

\begin{figure}
 \begin{center}
    \includegraphics[width=5cm]{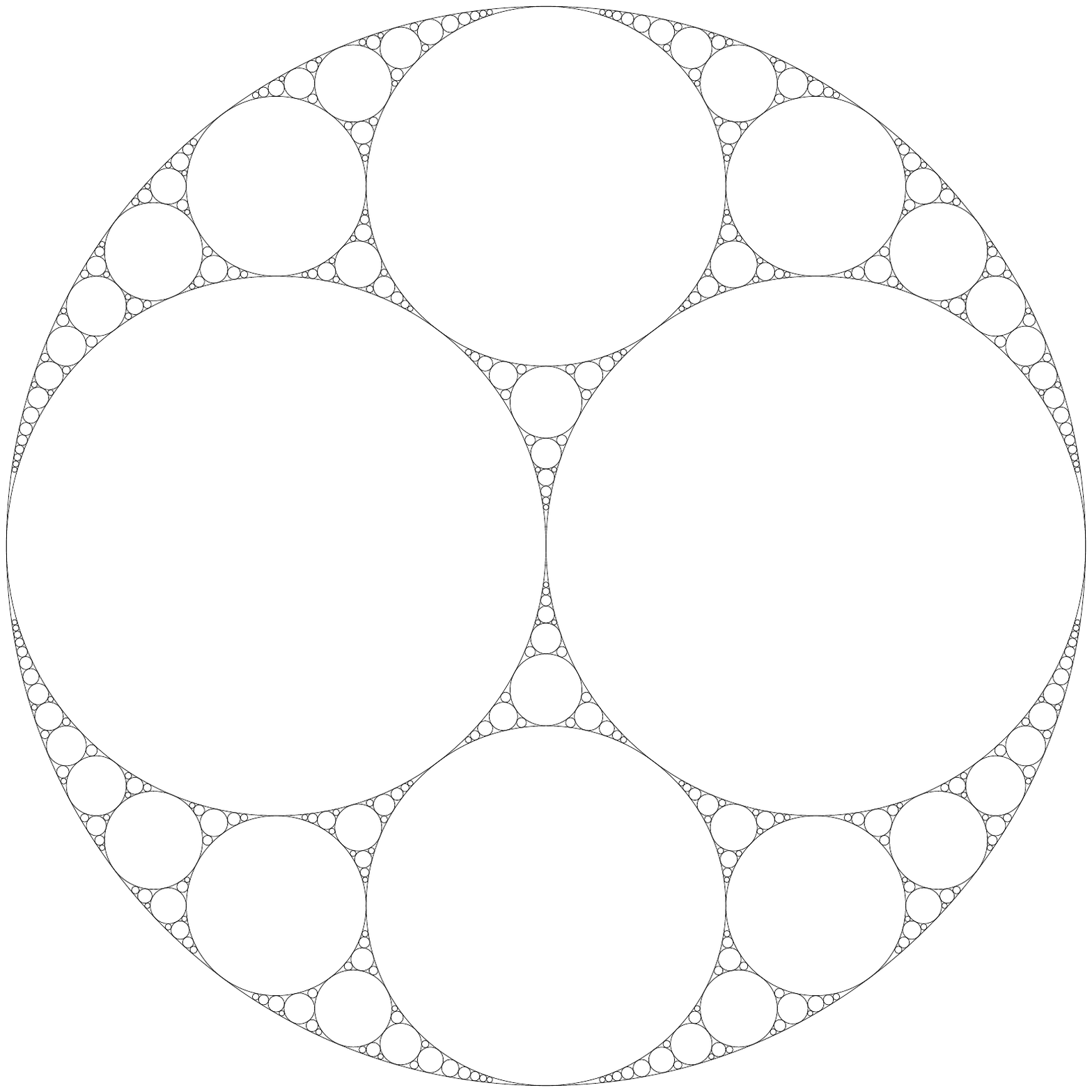}
    \includegraphics[width=5cm]{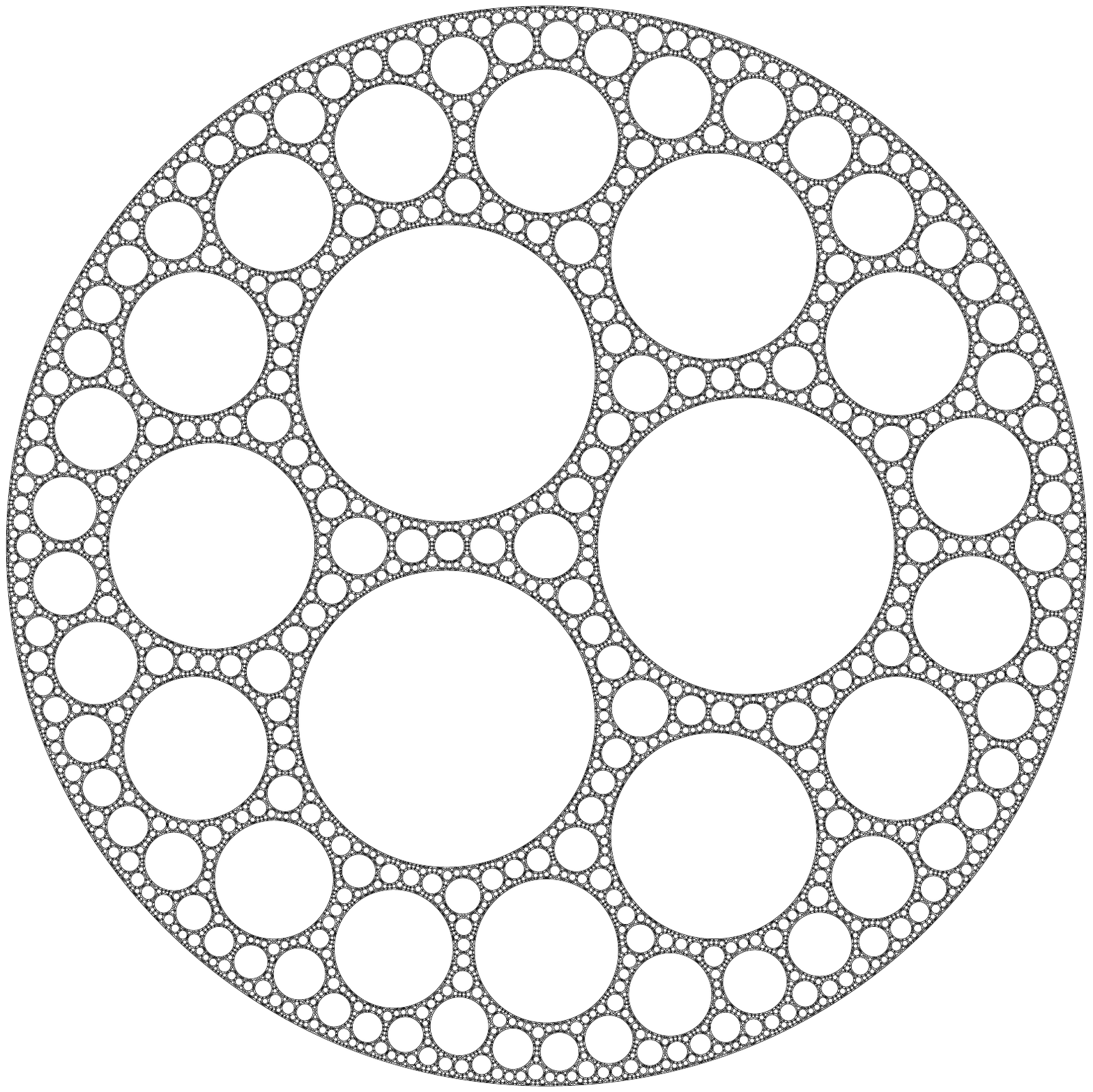}
    \caption{Apollonian circle packing and Sierpinski curve (by C. McMullen)}
    \label{f3}
\end{center}
\end{figure}

We suppose that $\mathcal P$ is locally finite in the sense that for any $T>1$, there
are only finitely many circles in $\P$ of curvature at most $T$ in any fixed bounded
subset of $\c$. Geometrically, $\P$ is locally finite if there
is no infinite sequence of circles in $\P$ converging to a fixed circle.
For instance, if the circles of $\P$ have disjoint interiors as in Fig.~\ref{f3}, $\P$ is locally finite.

For a bounded subset $E$ of $\c$ and $T>1$, we set
$$N_T(\P, E):=\# \{C\in \P: C\cap
E\ne\emptyset ,\;\;\op{Curv}(C)<T\}  $$
where $\op{Curv}(C)$ denotes the curvature of a circle $C$.
The local finiteness assumption on $\P$ implies that $N_T(\P, E)<\infty$.
Our question is then if there exists  a Borel measure $\omega_\P$ on $\c$
such that for all nice Borel subsets $E_1, E_2\subset \c$,
$$\frac{N_T(\P, E_1)}{N_T(\P, E_2)}\sim_{T\to \infty} \frac{\omega_\P(E_1)}{\omega_\P(E_2)},$$
assuming $N_T(\P, E_2)>0$ and $\omega_\P(E_2)>0$.

 Our main theorem applies to a very general packing
$\P$, provided $\P$ is invariant under a non-elementary (i.e., non virtually-abelian)
Kleinian group satisfying certain finiteness conditions.

Recall that a Kleinian group is a discrete subgroup of $G:=\PSL_2(\c)$
and $G$ acts on the extended complex plane $\hat \c=\c\cup\{\infty\}$ by M\"obius transformations:
$$\begin{pmatrix} a & b\\ c& d\end{pmatrix} z=\frac{az+b}{cz+d} $$
where $a,b,c,d\in \c$ with $ad-bc=1$ and $z\in \hat \c$.
A M\"obius transformation maps a circle to a circle and by the Poincare extension,
 $G$ can be identified with the group of all orientation preserving isometries of $\bH^3$.
Considering the upper-half space model $\bH^3=\{(z, r): z\in \c, r>0\}$,
the geometric boundary $\partial_\infty(\bH^3)$ is naturally identified with
$\hat \c$.

For a Kleinian group $\G$,
we denote by $\Lambda(\G)\subset\hat \c$ its limit set,
 that is, the set
of accumulation points of an orbit of $\G$ in $\hat \c$, and
 by $0\le \delta_\G\le 2 $ its critical exponent. For $\G$ non-elementary, it is known that
$\delta_\G >0$.
Let $\{\nu_x:x\in \bH^3\}$ be a $\G$-invariant conformal density  of
dimension $\delta_\G$ on $\Lambda(\G)$, which exists by the work
 of Patterson \cite{Patterson1976} and Sullivan \cite{Sullivan1979}.

In order to present our main theorem on the asymptotic of $N_T(\P, E)$
we introduce two invariants associated to $\G$ and $\P$.
The first one is a Borel measure on $\c$ depending only on $\G$.
\begin{Def}\rm  Define
a Borel measure $\omega_\G$ on $\c$: for $\psi\in C_c(\c)$
$$\omega_\G(\psi)=\int_{z\in \c} \psi (z) e^{\delta_\G \beta_z(x,(z,1)) }\; d\nu_{x}(z) $$
 where $x\in \bH^3$ and
$\beta_z(x_1,x_2)$ is the signed distance between the horospheres based at $z\in  \c$ and passing
through $x_1, x_2\in \bH^3$.

 By the conformal property of $\{\nu_{x}\}$,
 $\omega_{\G}$ is well-defined independent of the choice of $x\in \bH^n$.
\end{Def}
We have a simple formula: for $j=(0,1)\in \bH^3$,
$$d\omega_{\G} =(|z|^2+1)^{\delta_\G} d\nu_{j} .$$

 For a vector $u$ in the unit tangent bundle $\T^1(\bH^3)$, denote by $u^{+}\in \hat \c$
 (resp. $u^-\in \hat \c$) the
 forward (resp. backward) end point of the geodesic determined by $u$.
On the contracting horosphere $H^-_\infty(j)\subset \T^1(\bH^3)$ consisting of
 upward unit normal vectors on the horizontal plane $\{(z,1):z\in\c\}$,
 the normal vector based at $(z,1)$ is mapped to $z$ via
 the map $u\mapsto u^-$.
  Under this correspondence, the measure $\omega_\G$ on $\c$ is equal to the density of
   the Burger-Roblin measure $\tilde m^{\BR}$ (see Def. \ref{BR})
  on $H^-_\infty(j)$.

The second invariant is a number in $[0,\infty]$ measuring a certain size of $\P$.
\begin{Def}[The $\G$-skinning size of $\P$] {\rm
For a circle packing $\P$ invariant under $\G$,
define $0\le \op{sk}_\G(\P)\le \infty$ as follows:
$$\op{sk}_\G(\P):=\sum_{i\in I} \int_{s\in \op{Stab}_{\G} (C_i^\dagger)\ba C_i^\dagger}  e^{\delta_\G
\beta_{s^+}(x,\pi(s))}d\nu_{x}(s^+)$$
where $x\in \bH^3$, $\pi: \T^1(\bH^3)\to \bH^3$ is the canonical projection,
$\{C_i:i\in I\}$ is a set of representatives of $\G$-orbits in $\P$,
$C_i^\dagger\subset \T^1(\bH^3)$ is the set of unit normal vectors to
 the convex hull $\hat C_i$ of $C_i$ and $\op{Stab}_{\G} (C_i^\dagger)$ denotes the set-wise
 stabilizer of $C_i^\dag$ in $\G$. 
 Again by the conformal property of $\{\nu_x\}$, the definition of
 $\op{sk}_{\G}(\P)$ is independent of the choice of $x$ and the choice of representatives $\{C_i\}$.
}\end{Def}

We remark that the value of $\op{sk}_\G(\P)$ can be zero or infinite in general
and we do not assume any condition on $\op{Stab}_{\G} (C_i^\dagger)$'s (they may be trivial).

We denote by
 $m^{\BMS}_\G$ the Bowen-Margulis-Sullivan measure on the unit tangent bundle
 $\T^1(\Gamma\ba \bH^3)$ associated to the density $\{\nu_x\}$
(Def. \ref{bms}).
 When $\G$ is geometrically finite, i.e.,
$\G$ admits a finite sided fundamental domain in $\bH^3$,
Sullivan showed that $|m^{\BMS}_\G|<\infty$ \cite{Sullivan1984}
and that $\delta_\G $ is equal to
 the Hausdorff dimension of the limit set $\Lambda(\G)$ \cite{Sullivan1979}.
A point in $\Lambda(\G)$ is called a parabolic fixed point of $\G$ if it is fixed
by a parabolic element of $\G$.

\begin{figure}
 \begin{center}
    \includegraphics[width=5cm]{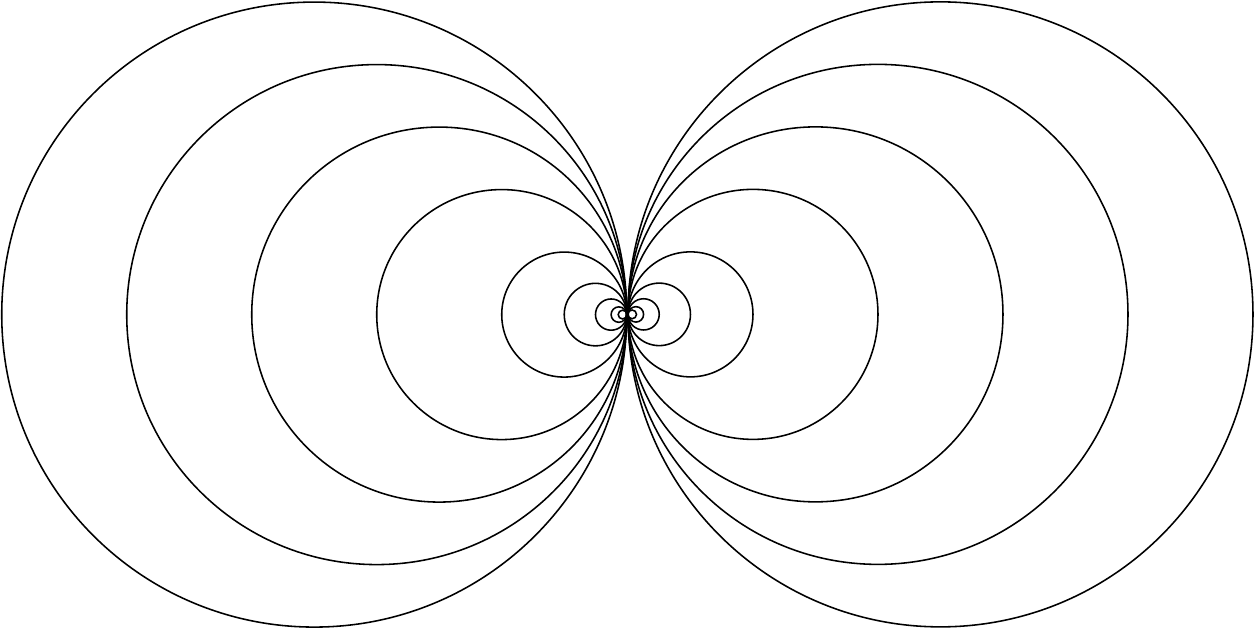}
    \caption{Infinite bouquet of tangent circles}
    \label{bouqqq}
\end{center}
\end{figure}

\begin{Def}\rm
 By an \emph{infinite bouquet of tangent circles glued at a point $\xi\in \c$},
 we mean a union of two collections, each consisting of 
infinitely many pairwise internally tangent circles with the common tangent point $\xi$ and their radii tending to $0$,
 such that the circles in each collection are externally tangent to the circles in the other at $\xi$ (see Fig.~\ref{bouqqq}).
\end{Def}

\begin{Thm}\label{m1gf}
 Let $\mathcal P$ be a
  locally finite circle packing in $\c$ invariant under
a non-elementary geometrically finite group $\G$ and with finitely many
$\G$-orbits.
If $\delta_\G\le 1$, we further assume that $\P$ does not contain an infinite
 bouquet of tangent circles glued at a parabolic fixed point of $\Gamma$.
Then $\op{sk}_{\G}(\P)<\infty$ and
for any bounded Borel subset $E$ of $\c$ with
$\omega_\G(\partial(E))=0$,
 $$\lim_{T\to \infty} \frac{N_T(\P, E)}{T^{\delta_\G}}=
 \frac{ \op{sk}_{\G}(\P) }{\delta_\G \cdot |m^{\BMS}_\G|} \cdot
\omega_\G(E)  . $$
If $\P$ has infinitely many circles, then $\op{sk}_{\G}(\P)>0$.
\end{Thm}

 \begin{Rmk}\rm
  \begin{enumerate}
\item Given a finite collection 
$\{C_1, \cdots, C_m\}$ of circles in the plane $\c$ and a non-elementary
 geometrically finite group $\G <\PSL_2(\c)$,
Theorem \ref{m1gf} applies to $\P:=\cup_{i=1}^m \G(C_i)$, provided
$\P$ contains neither infinitely many circles converging to a fixed circle 
nor any infinite bouquet of tangent circles.
\item In the case when $\delta_\G\le 1$ and
$\P$ contains an infinite bouquet of tangent circles glued at a parabolic fixed point of $\Gamma$, we have
$\op{sk}_{\G}(\P)=\infty$ \cite{OhShahGFH}.
In that case if the interior of $E$ intersects $\Lambda(\G)$
non-trivially, the growth order of $N_T(\mathcal P, E)$
is $T\log T$ if $\delta_\G=1$, and it is
$T$ if $\delta_\G<1$ \cite{OS}.
   \item We note that the asymptotic of $N_T(\P, E)$
depends only on $\G$, except for the $\G$-skinning size of $\P$.
This is rather surprising in view of the fact
that there are circle packings with completely different
configurations but invariant under the same group $\G$.
\item Theorem \ref{m1gf} implies that the asymptotic
distribution of small circles in $\P$ is completely determined by
the measure $\omega_\G$: for any bounded Borel
sets $E_1, E_2$ with $\omega_\G(E_2)>0$ and
$\omega_\G(\partial(E_i))=0$, $i=1,2$, as $T\to \infty$,
  $$\frac{N_T(\P, E_1)}{N_T(\P, E_2)} \sim \frac{\omega_\G(E_1)}{\omega_\G(E_2)} .$$

\item Suppose that all circles in $\P$ can be oriented so that they have
 disjoint interiors whose union is equal to 
 the domain of discontinuity $\Omega(\G):=\hat \c-\Lambda(\G)$. If
either $\P$ is bounded or $\infty$ is a parabolic fixed point
for $\G$, then
 $\delta_\G$ is equal to the circle packing exponent $e_{\mathcal P}$ defined as:
 $$e_{\mathcal P}=\inf\{s: \sum_{C\in \P} \op{Curv}^{-s}<\infty\}=\sup\{s:\sum_{C\in \P}\op{Curv}(C)^{-s}=\infty\}. $$
 
 This was proved by Parker~\cite{Parker1995} extending the earlier works of Boyd \cite{Boyd1973} and Sullivan \cite{Sullivan1984} on bounded
 Apollonian circle packings.
\end{enumerate}
 \end{Rmk}

In the proof of Theorem \ref{m1gf}, the geometric finiteness assumption on $\G$ is used only to
ensure the finiteness of the Bowen-Margulis-Sullivan measure $m^{\BMS}_\G$.
We have the following more general theorem:
\begin{Thm}\label{m1}
 Let $\mathcal P$ be a locally finite
   circle packing invariant under
 a non-elementary Kleinian group $\G$ and with finitely many $\G$-orbits.
Suppose that
$$|m_{\G}^{\BMS}|<\infty\quad\text{and}\quad  \op{sk}_\G(\P)<\infty .$$
Then for any bounded Borel subset $E$ of $\c$ with
$\omega_\G(\partial(E))=0$,
 $$\lim_{T\to \infty}\frac{N_T(\P, E)}{T^{\delta_\G}}
=\frac{ \op{sk}_{\G}(\P) }{\delta_\G \cdot |m^{\BMS}_\G|} \cdot
\omega_\G(E)  . $$
If $\P$ is infinite, then $\op{sk}_{\G}(\P)>0$.
\end{Thm}

Since there is a large class of geometrically infinite
groups with $|m^{\BMS}_\G|<\infty$ \cite{Peigne2003}, Theorem \ref{m1} is not subsumed
by Theorem \ref{m1gf}.

We remark that the condition on the finiteness of
$m_\G^{\BMS}$ implies that the density $\{\nu_x\}$ is determined uniquely up to homothety
(see \cite[Coro. 1.8]{Roblin2003}).

\begin{Rmk}{\rm
\begin{enumerate}
\item
The assumption of $|m_\G^{\BMS}|<\infty$ implies that
 $\nu_x$ (and hence $\omega_\G$) is atom-free \cite[Sec. 1.5]{Roblin2003}, and hence
the above theorem works for any bounded Borel subset $E$
intersecting $\Lambda(\G)$ only at finitely many points.

\item It is not hard to show that
$\G$ is Zariski dense in $\PSL_2(\c)$ considered as a real algebraic group if and only
if $\Lambda(\G)$ is not contained in a circle in $\hat \c$.
In such  a case,
any proper real subvariety of $\hat \c$ has zero $\nu_x$-measure.
This is shown in \cite[Cor.1.4]{FlaminioSpatzier} for $\G$ geometrically finite but
its proof works equally well if $\nu_x$ is $\G$-ergodic, which is the case when $|m_\G^{\BMS}|<\infty$.
Hence Theorem~\ref{m1} holds for any Borel subset $E$ whose boundary is contained in a countable union of real algebraic curves.
\end{enumerate}} \end{Rmk}

We now describe  some concrete applications of Theorem~\ref{m1gf}. 

\bigskip

\subsection{Apollonian gasket} 
\begin{figure}
 \begin{center}
    \includegraphics[width=4cm]{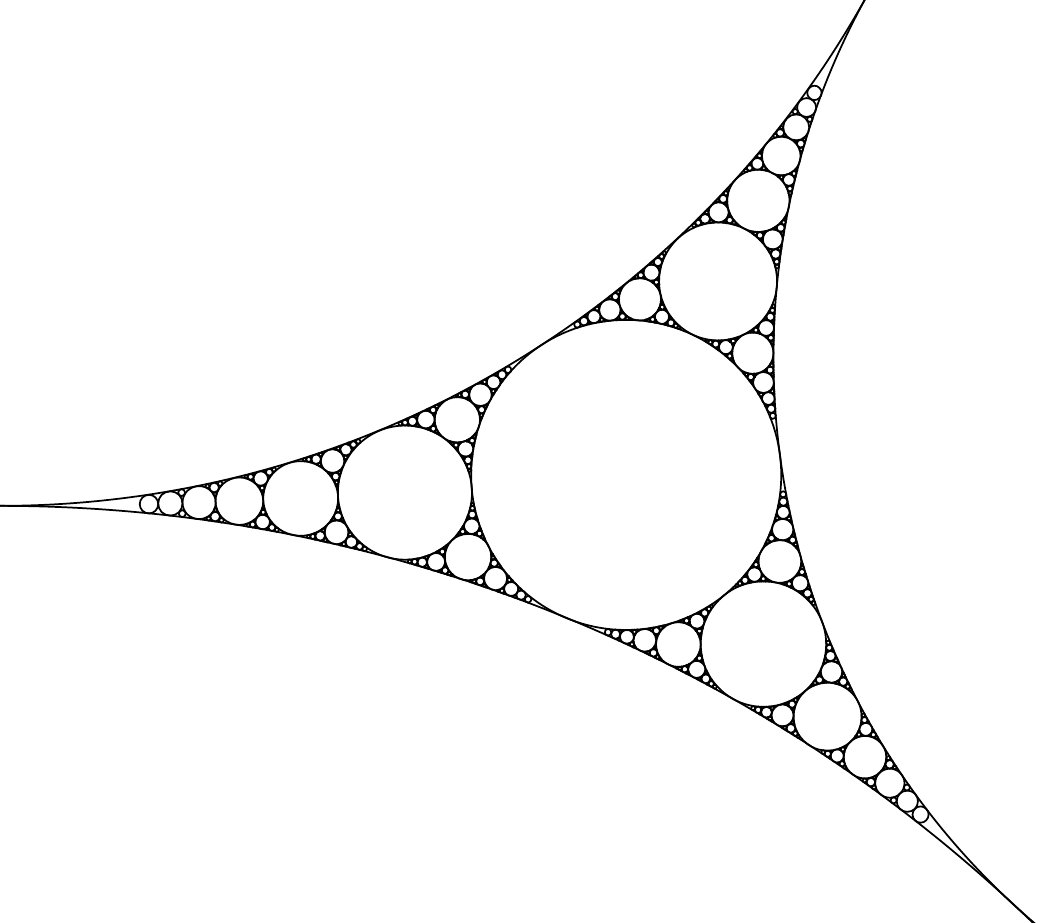}
     \caption{Apollonian gasket}\label{f4}
 \end{center}
\end{figure}

Three mutually tangent circles in the plane determine a curvilinear triangle, say, $\mathcal T$.
By a theorem of Apollonius of Perga (c. 200 BC), one can inscribe
a unique circle into the triangle $\mathcal T$, tangent to all of the
three circles. This produces three more curvilinear triangles inside $\mathcal T$ and
we inscribe a unique circle into each triangle. By continuing to add circles in this way,
we obtain an infinite circle packing of $\mathcal T$, called the Apollonian gasket for $\mathcal T$,
say, $\A$ (see Fig.~\ref{f4}).

By adding {\it all} the circles tangent to three of the given ones, not only those within $\mathcal T$,
one obtains an Apollonian circle packing $\mathcal P:=\mathcal P(\mathcal T)$, which may be bounded or unbounded (cf. \cite{GrahamLagariasMallowsWilksYanI} \cite{GrahamLagariasMallowsWilksYanI-n}, \cite{SarnakMAA}, \cite{SarnakToLagarias}, \cite{KontorovichOh}).

\begin{figure}\begin{center}
 \includegraphics [width=2in]{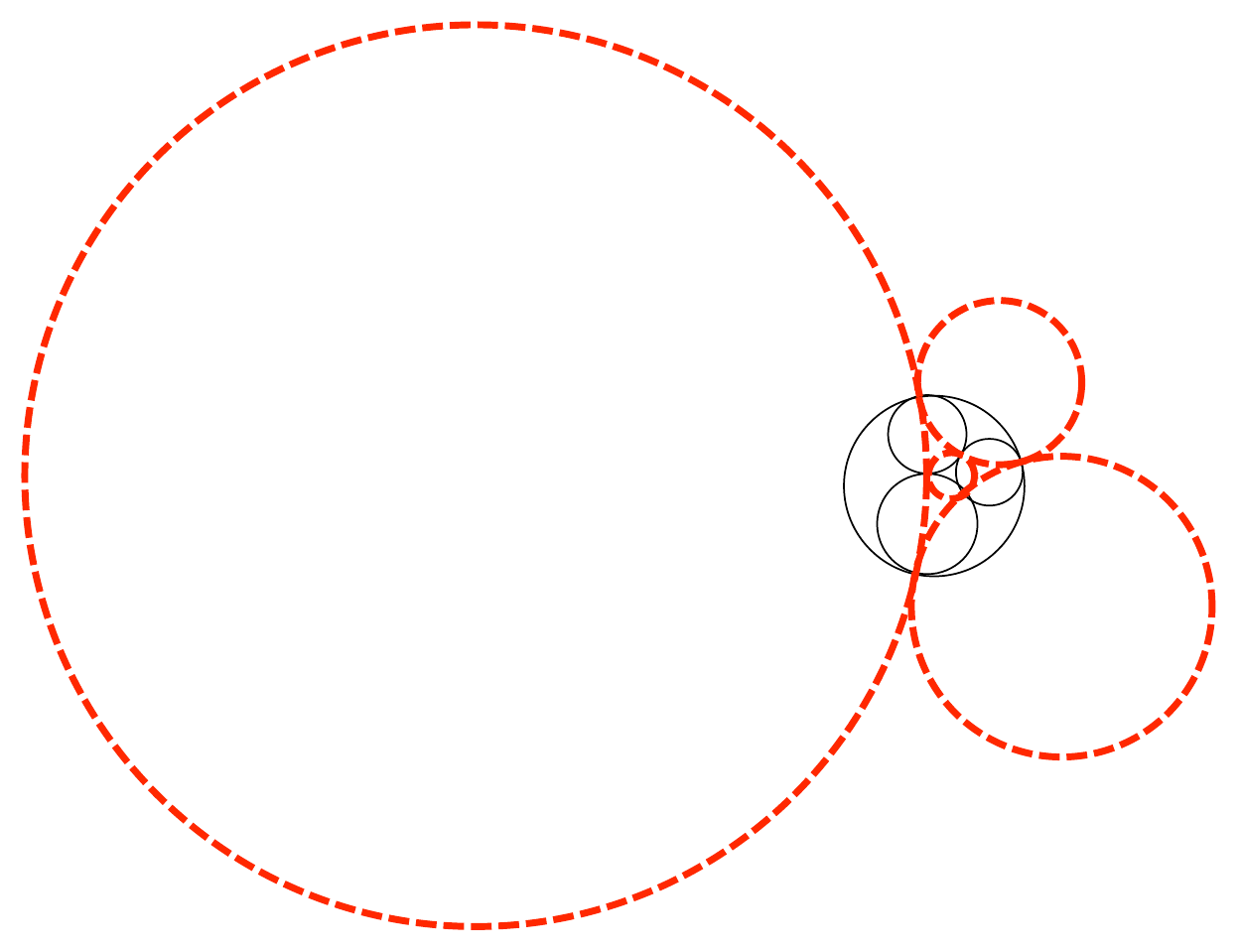}
\caption{Dual circles}
\label{picDual}\end{center}
\end{figure}

 Fixing four mutually tangent circles in $\P$,
consider the four dual circles determined by the six intersection points (see Fig. \ref{picDual}
where the dotted circles are dual circles to the solid ones),
 and
denote by $\G_\P$ the intersection of  $\PSL_2(\c)$ and the group
generated by the inversions with respect to those dual circles.
Then $\G_\P$ is a geometrically finite Zariski dense subgroup
of the real algebraic group $\PSL_2(\c)$ preserving $\P$, and
its limit set in $\hat \c$ coincides with the residual set of $\P$ (cf. \cite{KontorovichOh}).

We denote by $\alpha$
the Hausdorff dimension of the residual set of $\P$, which
is known to be $1.3056(8)$ according to McMullen \cite{McMullen1998}.

\begin{Cor}\label{Apo} Let $\mathcal{T}$ be a curvilinear triangle
determined by three mutually tangent circles and $\A$  the Apollonian gasket for
 $\mathcal{T}$.
Then for any Borel subset $E\subset \mathcal T$
whose boundary is contained in a countable union of real algebraic curves,
$$\lim_{T\to \infty}\frac{N_T(E)}{ T^{\alpha} }
 = \frac{ \op{sk}_{\G_{\P}}(\P) }{\alpha \cdot |m^{\BMS}_{\G_{\P}}|} \cdot
\omega_{\G_{\P}}(E)$$
 where $N_T(E):=\#\{C\in \A: C\cap E\ne \emptyset, \;\;  \op{Curv}(C)<T\}$ and
$\P=\P(\mathcal T)$.
\end{Cor}

Either when $\P$ is bounded
and $E$ is the disk enclosed by the largest circle of $\P$, or when $\P$ lies between two parallel lines
and $E$ is the whole period,
it was proved in \cite{KontorovichOh} that
$N_T(\P, E)\sim c  \cdot T^\alpha$ for some $c>0$. This implies that
{\footnote{$\asymp$ means that the ratio of the two sides is between two uniform constants}}
$N_T(\mathcal T)\asymp T^\alpha$.
The approach  in \cite{KontorovichOh}  was based on
the Descartes circle theorem in parameterizing quadruples of circles of curvature at most $T$ as
 vectors of maximum norm at most $T$
 in the cone  defined by the Descartes quadratic equation.
We remark that the fact that $\alpha$ is strictly bigger than $1$
was crucial in the proof of \cite{KontorovichOh} as based on the $L^2$-spectral theory
of $\G_{\mathcal P}\ba \bH^3$.

\begin{figure} 
\begin{center}
 \includegraphics [height=9cm]{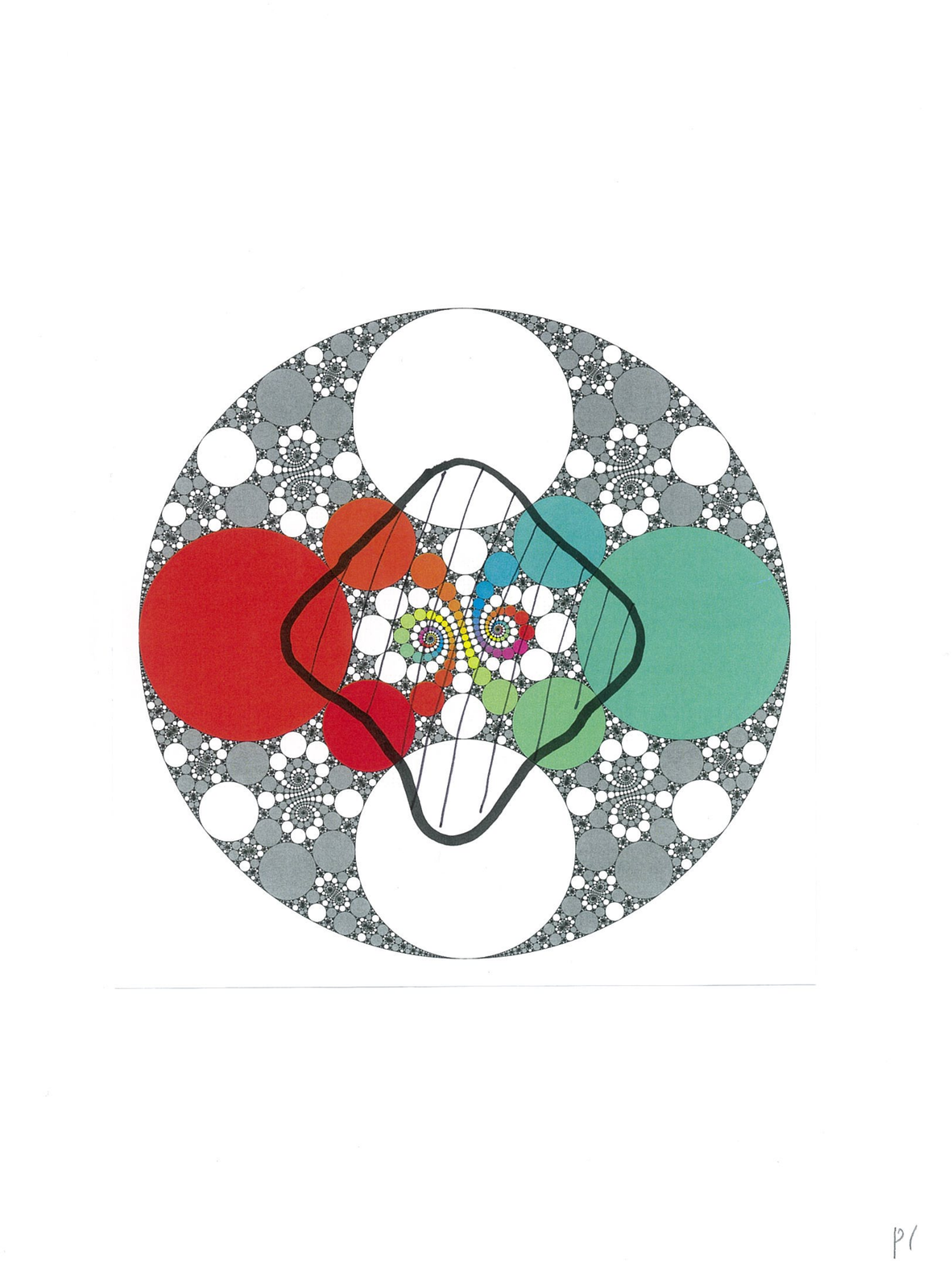}
 \includegraphics [height=8cm]{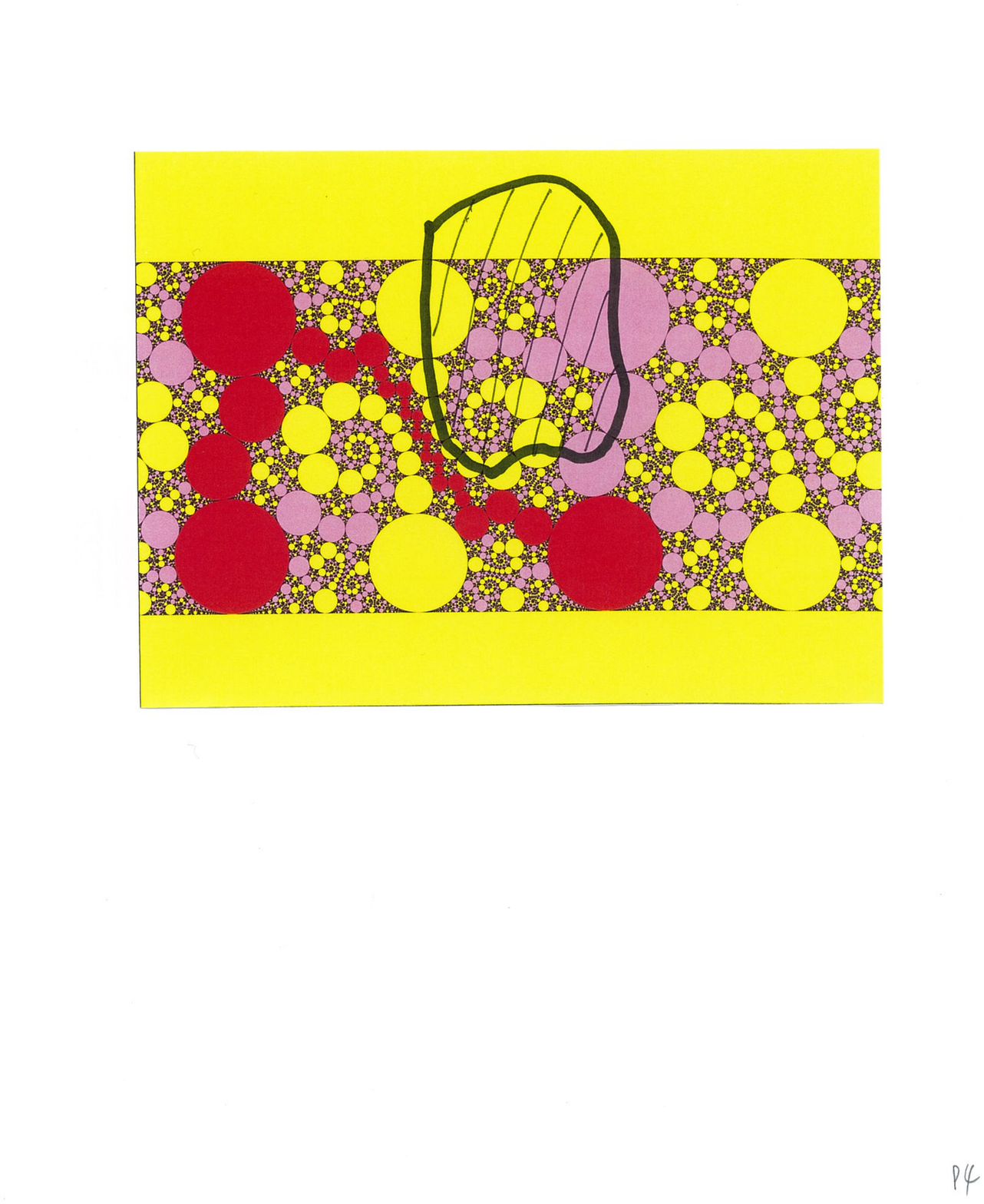}\end{center}
 \caption{Regions whose boundary intersects $\Lambda(\G)$ at finitely many points (background pictures are reproduced with permission from Indra's Pearls, by
D.Mumford, C. Series and D. Wright, copyright Cambridge University Press 2002)}
    \label{f2}
\end{figure}

\bigskip

\subsection{Counting circles in the limit set $\Lambda(\G)$}
If $X$ is a finite volume hyperbolic $3$-manifold with totally geodesic boundary, then
its fundamental group $\G:=\pi_1(X)$ is geometrically finite and $X$ is homeomorphic to
$\Gamma \ba \bH^3\cup \Omega(\G)$ where
 $\Omega(\G):= \hat \c -\Lambda(\G)$ is the domain
 of discontinuity \cite{Kojima1992}.
The set $\Omega(\G)$
 is a union of countably many disjoint open disks in this case and has finitely many
 $\G$-orbits by the Ahlfors finiteness theorem \cite{Ahlfors1964}.
  Hence Theorem \ref{m1gf} applies to counting
 these open disks in $\Omega(\G)$ with respect to the curvature.
 
For example, for the group $\G$ generated by reflections in the sides of a unique
 regular tetrahedron whose convex core is bounded by four $\frac{\pi}4$ triangles and
 four right hexagons, $\Omega(\G)$ is illustrated 
in the second picture in Fig. \ref{f3} (see \cite[P.9]{McMullennotes275} for details).
 This circle packing is called a {\it Sierpinski curve}, being homeomorphic to
 the well-known Sierpinski carpet \cite{Claytor1934}.

Two pictures in Fig. \ref{f2} can be found
 in the beautiful book {\it Indra's pearls} by Mumford, Series and Wright
 (see P. 269 and P. 297
of \cite{MumfordSeriesWright})
where one can find many more circle packings to which our theorem applies.
 The book presents explicit geometrically finite Schottky groups $\G$ whose limit sets
 are illustrated in Fig. \ref{f2}.
The boundaries of the shaded regions meet $\Lambda(\G)$ only at finitely many points.
Hence our theorem applies to counting circles in these shaded regions.

\begin{figure}
 \begin{center}
    \includegraphics[width=6cm]{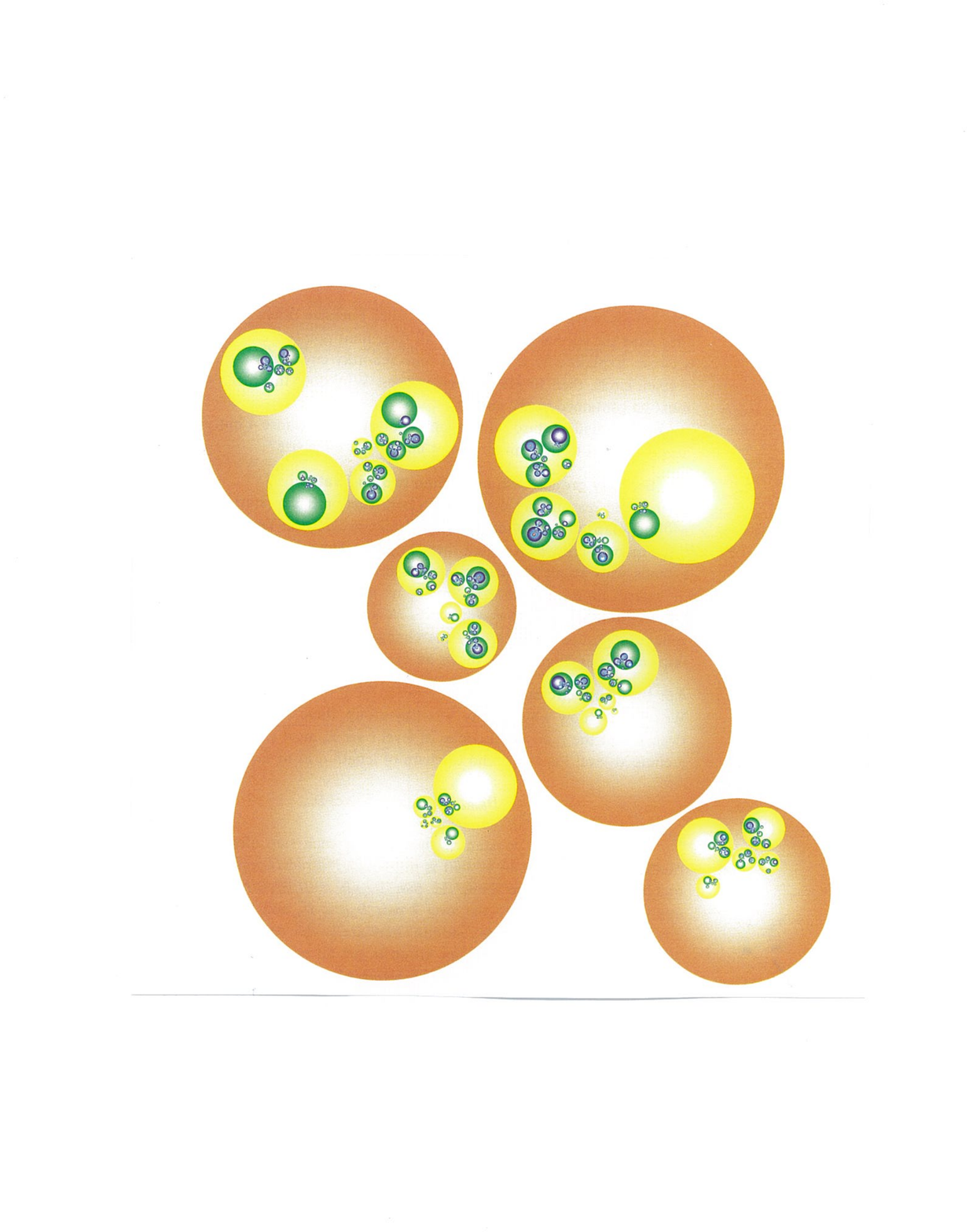}
     \caption{Schottky dance (reproduced with permission from Indra's Pearls, by
D.Mumford, C. Series and D. Wright, copyright Cambridge University Press 2002)} \label{f7}
 \end{center}
\end{figure}

\subsection{Schottky dance}
Another class of examples is obtained by considering the images of
Schottky disks under Schottky groups.
Take $k\ge 1$ pairs of mutually disjoint closed disks
$\{(D_i, D_i') : 1\le i\le k\}$ in $\c$ and for each $1\le i\le k$, choose
 a M\"obius transformation $\gamma_i$ which maps the interior of
$D_i$ to the exterior of $D_i'$ and the interior of $D_i'$ to
the exterior of $D_i$. The group, say, $\G$,
 generated by $\{\gamma_i:1\le i\le k\}$ is called a Schottky group of genus $k$ (cf. \cite[Sec. 2.7]{MardenOutercircles}).
The $\G$-orbit of the disks $D_i$ and $D_i'$'s nests
 down onto the limit set $\Lambda(\G)$ which is totally disconnected.
If we denote by $\P$ the union $\cup_{1\le i\le k} \G(C_i)\cup \G(C_i')$
where $C_i$ and $C_i'$ are the boundaries of $D_i$ and $D_i'$ respectively,
then $\P$ is locally finite, as the nesting disks become smaller and smaller.
The common exterior of hemispheres above the initial disks $D_i$ and $D_i'$
is a fundamental domain for $\G$ in the upper half-space $\bH^3$, and hence
$\G$ is geometrically finite.
Since $\P$ contains no infinite bouquet of tangent circles,
Theorem \ref{m1gf} applies to $\P$; for instance,
we can count circles in the picture in Fig. \ref{f7} (\cite[Fig. 4.11]{MumfordSeriesWright}).

\subsection*{On the structure of the proof}
In \cite{KontorovichOh}, the counting problem for a bounded Apollonian circle packing
was related to the equidistribution of expanding closed horospheres
on the hyperbolic $3$-manifold $\G\ba \bH^3$. For a general circle packing,
there is no analogue of the Descartes circle theorem which made such a relation possible.
The main idea in our paper
is to relate the counting problem for a general circle packing $\P$ invariant
under $\G$ with 
the equidistribution of orthogonal translates of a closed totally geodesic surface
in $\T^1(\G\ba \bH^3)$.
Let $C_0$ denote the unit circle centered at the origin and $H$ the stabilizer of $C_0$ in $\PSL_2(\c)$.
 Thus $H\ba G$ may be considered as the space of totally geodesic planes of $\bH^3$.
 The important starting point is to describe certain subset $B_T(E)$ in 
 $H\ba G$ so that
  the number of circles in the packing $\P:=\G(C_0)$ of curvature at most $T$ intersecting $E$
  can be interpreted as the number of
 points in $B_T(E)$ of a discrete $\G$-orbit on $H\ba G$.
  We then describe the weighted limiting distribution of orthogonal
  translates of an $H$-period $(H\cap \G)\ba H$ (which
   corresponds to a properly immersed hyperbolic surface which may be of infinite area)
    along these sets $B_T(E)$ in terms of the Burger-Roblin measure (Theorem \ref{mtt}) using the main result in \cite{OhShahGFH}
    (see Thm. \ref{os}). 
To translate the weighted limiting distribution result into the asymptotic
for $N_T(\P, E)$, we
relate the density of the Burger-Roblin measure of the contracting horosphere $H_\infty^-(j)$
  with the measure $\omega_\G$.

A version of Theorem~\ref{m1gf} in a weaker form, and some of its applications
 stated above were announced in \cite{OhICM}. We remark that the
methods of this paper can be easily generalized to prove a similar result
for a sphere packing in the $n$-dimensional Euclidean space invariant
under a non-elementary discrete subgroup of $\op{Isom}(\bH^{n+1})$.
\subsection*{Acknowledgment} We would like to thank Curt McMullen
for inspiring discussions. The applicability of our other paper \cite{OhShahGFH} 
in the question addressed in this paper came
 up in the conversation
of the first named author with him during her one month visit to Harvard in October, 2009.
She thanks the Harvard mathematics department for the hospitality.
We would also like to thank Yves Benoist for helpful discussions.

\section{Expansion of a hyperbolic surface by orthogonal geodesic flow}
\label{sectionh}
We use the following coordinates for the upper half space model for $\bH^3$:
$$\bH^3=\{z+rj=(z,r) :z\in \c, r>0\}$$
where $j=(0,1)$.
The isometric action of $G=\PSL_2(\c)$, via the Poincare extension of
the linear fractional transformations, is explicitly given as the following (cf. \cite{ElstrodtGrunewaldMennickebook}):
 \begin{equation}\label{EGM}
\begin{pmatrix} a & b\\ c& d\end{pmatrix} (z+rj)=\frac{(az+b)(\bar c\bar z +\bar d)+a\bar c r^2}{|cz+d|^2+|c|^2r^2}
 +\frac{r}{|cz+d|^2+|c|^2r^2} \;j .\end{equation}
 In particular, the stabilizer of $j$ is the following maximal compact subgroup of $G$:
 $$K:=\op{PSU}(2)=\{\begin{pmatrix} a & b\\ -\bar b & \bar a \end{pmatrix}: |a|^2+|b|^2=1 \} .$$

We set
$$ A:=\{a_t:=\begin{pmatrix} e^{t/2} & 0\\ 0& e^{-t/2}\end{pmatrix}:t\in \br\},\quad M:=\{ \begin{pmatrix} e^{i\theta} & 0\\ 0& e^{-i\theta}\end{pmatrix} :\theta\in \br \} $$
and $$
N:=\{n_z :=\begin{pmatrix} 1 & z\\ 0& 1\end{pmatrix}:z\in \c\} ,\quad
N^-=\{n_{z}^-:=\begin{pmatrix} 1 & 0\\ z& 1\end{pmatrix}:z\in \c\} .$$

We can identify $\bH^3$ with $G/K$  via the map $g(j)\mapsto gK$.
 Denoting by $X_0\in \T^1(\bH^3)$ the upward unit normal vector based at $j$,
 we can also identify the unit tangent bundle $\op{T^1}(\bH^3)$ with $G.X_0=G/M$:
here $g.X_0$ is given by $d\lambda(g) (X_0)$ where $\lambda(g): G/K\to G/K$ is the left translation $\lambda(g)(g'K)=
gg'K$ and $d\lambda(g)$ is its derivative at $j$.

 The geodesic flow $\{g^t\}$ on
 $\T^1(\bH^3)$ corresponds to the right translation by $a_t$ on $G/M$:
 $$g^t(gM)=ga_tM .$$

For a circle $C$ in $\c$, denote by $\hat C$  its convex hull, which is the northern hemisphere above $C$.

Set $C_0$ to be the unit circle in $\c$ centered at the origin.
 The set-wise stabilizer of $\hat C_0$ in $G$ is given by
 $$H=\op{PSU}(1,1) \cup \begin{pmatrix} 0&1 \\ -1 & 0\end{pmatrix}\op{PSU}(1,1)$$
 where
 $$\op{PSU}
 (1,1)=\{ \begin{pmatrix} a&b \\ \bar b &\bar a\end{pmatrix} : |a|^2-|b|^2=1 \} .$$
Note that $H$ is equal to the stabilizer of $C_0$ in $G$ and hence
 $\hat C_0$ can be identified with $H/H\cap K$.

We have the following generalized Cartan decomposition (cf. \cite{Schlichtkrull1984}):
for $A^+=\{a_t: t\ge 0\}$,
$$G=HA^+K $$
in the sense that every element of $g\in G$ can be written as $g=hak$, $h\in H, a\in A^+, k\in K$
and $h_1a_1k_1=h_2a_2k_2$ implies that $a_1=a_2$, $h_1= h_2 m$ and $k_1= m^{-1} k_2$ for
some $m\in H\cap K\cap Z_G(A)=M$.

As $X_0$ is orthogonal to the tangent space $\T_{j}(\hat C_0)$,
 $H.X_0=H/M$ corresponds to the set of unit normal vectors
   to $\hat C_0$, which we will denote by $C_0^\dagger$. Note that $C_0^\dagger$ 
has two connected
components, depending on their directions. 
For $t\in \br$, the set $g^t( C_0^\dagger)=(H/M)a_t=(Ha_tM)/M$
 corresponds to a union of two surfaces consisting of the orthogonal translates of $\hat C_0$ by distance $|t|$ in each direction, both having the same boundary $C_0$. 
\bigskip

Let $\G<G$ be a non-elementary discrete subgroup.
As in the introduction, let
 $\{\nu_x:x\in \bH^3 \}$ be a $\G$-invariant conformal density on $\hat \c$ of dimension $\delta_\G$,
  that is, each $\nu_x$ is a finite measure on $\hat \c$ 
  satisfying that
  for any $x,y\in \bH^3$, $z\in \hat \c$ and $\gamma\in
\G$,
$$\gamma_*\nu_x=\nu_{\gamma x};\quad\text{and}\quad
 \frac{d\nu_y}{d\nu_x}(z)=e^{-\delta_\G \beta_{z} (y,x)}. $$
Here $\gamma_*\nu_x(R)=\nu_x(\gamma^{-1}(R))$ for a Borel subset $R\subset \hat{\c}$
and the Busemann function $\beta_z(y_1, y_2)$ is given by
$\lim_{t\to\infty} d(y_1, \xi_t)-d(y_2,\xi_t)$ for a geodesic ray
$\xi_t$ toward  $z$.

 For $u\in \T^1(\bH^3)$, we define $u^{+}\in \hat \c$
 (resp. $u^-\in \hat \c$) to be the
 forward (resp. backward) end point of the geodesic determined by $u$ and $\pi(u)\in \bH^3$ to be the basepoint. Fixing $o\in \bH^3$, the map $u\mapsto (u^+, u^-,
t:=\beta_{u^-}(\pi(u), o))$ is
a homeomorphism between $\T^1(\bH^3)$ and $(\hat \c\times \hat \c - \{(\xi,\xi):\xi\in
\hat \c\}) \times \br$.
\begin{Def}\label{bms} \rm The Bowen-Margulis-Sullivan measure $m^{\BMS}_\G$ associated to $\{\nu_x\}$
(\cite{Bowen1971}, \cite{Margulisthesis}, \cite{Sullivan1984}) is the measure on $\T^1(\G\ba \bH^3)$
 induced by the following
$\G$-invariant measure on $\T^1(\bH^3)$: for $x\in \bH^3$,
$$d  \tilde m^{\BMS}(u)=e^{\delta_\G \beta_{u^+}(x, \pi(u))}\;
 e^{\delta_\G \beta_{u^-}(x,\pi(u)) }\;
d\nu_{x}(u^+) d\nu_{x}(u^-) dt .$$
\end{Def}
By the conformal properties of $\{\nu_x\}$,
 this definition is independent of the choice of $x\in \bH^3$.

We also denote by $\{m_x:x\in \bH^3\}$ a
$G$-invariant conformal
density of dimension $2$, which is unique up to homothety:
 each $m_x$ a finite measure on $\hat\c$ which is invariant under
$\op{Stab}_G(x)$ and $d m_x(z) =e^{-2 \beta_z(y,x)} dm_y(z)$ for any $x,y\in \bH^3$ and $z\in \hat \c$.

\begin{Def}\label{BR} {\rm  The Burger-Roblin measure $m^{\BR}_\G$
associated to $\{\nu_x\}$ and $\{m_x\}$ (\cite{Burger1990}, \cite{Roblin2003})
is the  measure on $\T^1(\G\ba \bH^n)$ induced by the following $\G$-invariant measure on $\T^1(\bH^n)$:
$$d  \tilde m^{\BR}(u)=e^{2\beta_{u^+}(x, \pi(u))}\;
 e^{\delta_\G \beta_{u^-}(x,\pi(u)) }\;
dm_x(u^+) d\nu_x(u^-) dt $$
for $x\in \bH^3$.
By the conformal properties of $\{\nu_x\}$ and $\{m_x\}$,
 this definition is independent of the choice of $x\in \bH^3$. 
} \end{Def}

For any circle $C$, let 
\[
H_C=\{g\in G:gC=C\}=\{g\in G: gC^\dag=C^\dag\}.
\]
We consider the following two measures on $C^\dag$: Fix any  $x\in \bH^3$, and let
\begin{equation}\label{hars}
d\mu^{\Leb}_{C^\dag} (s):=e^{2\beta_{s^+}(x,\pi(s))}dm_x(s) \text{ and }
 d\mu^{\PS}_{C^\dag} (s):=e^{\delta_\G \beta_{s^+}(x,\pi(s))} d\nu_x(s^+) .
\end{equation}
These definitions are independent of the choice of $x$ and
$\mu^{\Leb}_{C^\dag}$ (resp. $\mu^{\PS}_{C^\dag}$) is left-invariant by $H_C$ (resp. $H_C\cap \G)$).
Hence we may consider the measures
$\mu^{\Leb}_{C^\dag}$ and $\mu^{\PS}_{C^\dag}$  on the quotient $(H\cap \G)\ba C^\dag$.

We denote by $\op{sk}_\G(C)$ the total mass of  $\mu^{\PS}_{C^\dag}$; that is,
$$\op{sk}_\G(C):=\int_{s\in \G\cap H\ba C_0^\dag} e^{\delta_\G \beta_{s^+}(x,\pi(s))} d\nu_x(s^+) .$$
In general, $\op{sk}_\G(C)$ may be zero or infinite.

\begin{Thm}[{\cite[Theorem~1.9]{OhShahGFH}}] \label{os}
Suppose
  that the natural projection map $\G\cap H_C \ba \hat C \to \G\ba \bH^3$ is proper.
Assume that
$|m_\G^{\BMS}|<\infty$ and 
$\op{sk}_\G(C)<\infty$. Then for any $\psi\in C_c(\G\ba G/M)$,
as $t\to  \infty$,
$$e^{(2-\delta_\G)t}\int_{s\in (\G\cap H_C)\ba C^\dagger}\psi (sa_t) d\mu^{\Leb}_{C^\dag} (s)
\sim\frac{\op{sk}_\G(C) }{|m^{\BMS}_\G|}
m^{\BR}_\G(\psi) . $$
Moreover $\op{sk}_\G(C) >0$ if $[\G: H_C\cap \G]=\infty$.
\end{Thm}

Note that if
$|m_\G^{\BMS}|<\infty$, then $\G$ is of divergence type; that is, the Poincare series
of $\G$ diverges at $\delta_\G$. When $\Gamma$ is of divergence type,
 the $\G$-invariant conformal density $\{\nu_x\}$
of dimension $\delta_\G$ is unique up to homothety
 (see \cite[Remark following Corollary 1.8]{Roblin2003}):
explicitly $\nu_x$ can be taken as the weak-limit as $s\to \delta_\G^+$
of the family of measures
$$\nu_{x}(s):=\frac{1}{\sum_{\gamma\in \G} e^{-sd(j, \gamma j)}}
\sum_{\gamma\in\G} e^{-sd(x, \gamma j)} \delta_{\gamma j} .$$

Recall that $g\in \PSL_2(\c)$ is parabolic 
if and only if $g$ has a unique fixed point in $\hat \c$. 

\begin{Thm}[{\cite[Theorem~5.2]{OhShahGFH}}] \label{bou} Let $\G$ be geometrically finite. 
Suppose
  that the natural projection map $\G\cap H_C \ba \hat C \to \G\ba \bH^3$ is proper.
Then $\op{sk}_\G(C)<\infty$ if and only if  either $\delta_\G > 1$ or
any parabolic fixed point of $\G$ lying on $C$ is
 fixed by a parabolic element of $H_C\cap \G$. 
\end{Thm}

\begin{proof}
Note that in the notation of  
{\cite[Theorem~5.2]{OhShahGFH}},
 if we put $E=\hat C$, which is a complete totally geodesic submanifold of $\bH^3$ of
codimension $1$, 
then $\partial(\pi(\tilde E))=C$, $\tilde E=C^\dag$, $\Gamma_{\tilde E}=H_C\cap \G$, 
and $\lvert \mu_E^{\PS}\rvert=\op{sk}_{\G}(C)$. Hence the conclusion is immediate.
\end{proof}

 \section{Reformulation into the orbital counting problem on the space of hyperbolic planes}
 Let $G=\PSL_2(\c)$ and
  $\G<G$ be a non-elementary discrete subgroup. Let $C$ be a circle in $\hat \c$
  and $H_C$ denote the set-wise stabilizer of $C$ in $G$.

It is clear:
\begin{Lem}\label{infinite}
If $\G(C)$ is infinite, then $[\G:H_C\cap \G]=\infty$.
\end{Lem}

\begin{Lem}\label{locfinite}\label{lc}
The following are equivalent:
\begin{enumerate}
 \item 
 A circle packing $\G(C)$ is locally finite;
\item
the natural projection map $f:\G\cap H_C \ba \hat C \to \G\ba \bH^3$ is proper;
\item
$H_C\ba H_C\G $ is discrete in $H_C\ba G$.
\end{enumerate}
\end{Lem}
\begin{proof}
 We observe that the properness of  $f$ is equivalent to the condition that
only finitely many distinct hemispheres in $\G(\hat C)$ intersects a given compact subset of $\bH^3$.
Note that any compact subset of $\bH^3$ is contained in a compact subset 
the form $E\times [r_1, r_2]=\{(z,r): z\in E, \, r_1\le r\le r_2\}$ for $E\subset \c$
 compact and $0< r_1 < r_2<\infty$, and that the radius of a circle in $\c$ is same as
the height of its convex hull in $\bH^3$. Hence the properness of the map $f$ is again
equivalent to the condition that for any $r>0$ and a compact subset $E\subset \c$,
 there are only finitely many distinct circles in $\G(C)$
intersecting $E$ and of radii at least $r$, that is, $\G(C)$ being locally finite,
proving the equivalence of (1) and (2).

It is straightforward to verify that the properness of $f$ and that of the projection
 map $\Gamma\cap H_C\ba C^\dag\to \Gamma\ba \T^1(\bH^3)$ are equivalent.
Let $X_C\in C^\dag$ such that $X_C^+=\infty\in \hat\c$. 
Let $M_C=\{g\in G: gX_C=X_C\}$. Since $\hat{C}$ is the unique totally geodesic 
submanifold of  $\bH^3$ orthogonal to $X_C$, $M_C$ is contained in $ H_C$. 
We identify $G/M_C$ with $\T^1(\bH^3)$ via $gM_C\mapsto gX_C$. Since $H/M_C$ 
identifies with $C^\dag$, the canonical map $\Gamma\cap H_C\ba  H_C/M_C\to \Gamma\ba G/M_C$ 
is proper. Since $M_C$ is compact, it follows
 that $\G\cap H_C\ba H_C\to \G\ba G$ is proper. Equivalently $\G H_C$ is closed in $G$ (see \cite{OhShahGFH} for the equivalence). 
As $\G$ is countable, 
this is again equivalent to the discreteness of $H_C\ba H_C\G$ in $H_C\ba G$. This proves the
equivalence of (2) and (3).
\end{proof}

\begin{Rmk}\rm  If $\G\cap H_C$ is a lattice in $H_C$, then $\G H_C$ is closed in $G$ (\cite[\S1]{msrbook}), and hence
$\G(C)$ is a locally finite circle packing. In this case, by \cite[Theorem~1.11]{OhShahGFH}, we have
$\op{sk}_\G(C)<\infty$. 
\end{Rmk}

\begin{Prop}\label{pf} Let $\xi\in C$ be a parabolic fixed point of $\G$.
Suppose that $\G(C)$ does not contain an infinite bouquet of tangent circles glued at $\xi$. 
Then $\xi$ is a parabolic fixed point for $H_C\cap \G$.
\end{Prop}

\begin{proof} Suppose that there exists a parabolic element $\gamma\in \G - H_C$ 
fixing $\xi\in\c$. By sending $\xi$ to $\infty\in\hat\c$ by an element of $G$, 
we may assume that $\xi=\infty$ and
$\gamma$ acts as a translation on $\c$.  
Since $\gamma C\neq C$ and $C$ is a circle passing through $\infty$, 
we have that $\{\gamma^{k}C:k\in\z\}$ is an infinite collection of parallel lines. 
By sending $\infty$ back to the original $\xi$, we see that $\{\gamma^{k}C:k\in\z\}$ is an infinite bouquet of tangent circles glued at $\xi$.
\end{proof}

\subsection{Deduction of Theorem~\ref{m1gf} from Theorem~\ref{m1}} 
We only need to ensure that $\sk_{\G}(\P)<\infty$, or equivalently, $\sk_\G(C)<\infty$ for each $C\in\P$. 
By the assumption in Theorem~\ref{m1gf}, if $\xi\in C$ is any parabolic fixed point of $\Gamma$, then by Proposition~\ref{pf}, $\xi$ is a parabolic fixed point for $H_C\cap \G$. Therefore by Theorem~\ref{bou}, $\op{sk}_\G(C)<\infty$. 
\qed

\subsection{Relating counting on a single $\G$-orbit to a set $B_T(E)\subset H\ba G$}

In the rest of this section,  let $C_0$ denote the unit circle in $\c$ centered
at the origin
and let $H:=\op{Stab}(\hat C_0)$. We follow notations from Section \ref{sectionh}.
We assume that
$\G(C_0)$ is a locally finite circle packing of $\c$.

Let $E$ be a bounded subset in $\c$
and set
$$N_T(\G(C_0), E):=\#\{C\in \G(C_0): C\cap E\ne \emptyset,\quad \op{Curv}(C)<T\} .$$

For $s>0$, we set
 $$A^+_s:=\{a_t: 0\le t\le s\};\quad  A^-_s:=\{a_{-t}: 0\le t\le s\} .$$
For a subset $E\subset \c$, we set $N_E:=\{n_z:z\in E\}$.
\begin{Def}[Definition of $B_T(E)$]
{\rm For $E\subset \c$ and $T>1$, we define the subset $B_T(E)$ of $H\ba G$ to be the image of
the set $$KA^+_{\log T}N_{-E}=\{ ka_t n_{-z}\in G:  k\in K, 0\le t< \log T, z\in E\}$$
under the canonical projection $G\to H\ba G$.}
\end{Def}

For a bounded circle $C$ in $\c$, $C^\circ$ denotes the open disk enclosed by $C$.
We will not need this definition for a line since
there can be only finitely many lines intersecting a fixed bounded
subset in a locally finite
circle packing.

\begin{Def}\label{adm}{\rm For a given circle packing $\P$,
 a bounded subset $E\subset \c$ is said to be \emph{$\mathcal P$-admissible}\/ if,
  for any bounded circle $C\in \mathcal P$, $C^\circ\cap E\ne \emptyset$ implies
  $C^\circ\subset E$, possibly except for finitely many circles.
}\end{Def}

The following translation of $N_T(\G (C_0), E)$ as
the number of points in $[e]\G \cap B_T(E)$, where $[e]=H\in H\ba G$,
 is crucial in our approach:

\begin{Prop}\label{tran} If $E$ is $\G(C_0)$-admissible,
there exists $m_0\in \mathbb N$ such that for all $T\gg 1$,
$$ \# [e]\G\cap B_T(E) - m_0 \le N_T(\G (C_0), E)\le \# [e]\G\cap B_T(E) + m_0.$$
\end{Prop}
\begin{proof}
Observe that
\begin{align*}
\# [e]\G\cap B_T(E)&=\#\{[\gamma]\in \G\cap H\ba \G: H\gamma\cap KA_{\log T}^+N_{-E}\ne\emptyset \} \\
&=\# \{[\gamma]\in \G/\G\cap H:   \gamma HK\cap N_E A_{\log T}^-K\ne \emptyset\}  \\ &=\#\{\gamma (\hat C_0): \gamma HK\cap N_EA_{\log T}^-K\ne \emptyset\} \end{align*}
where the second equality is obtained by taking the inverse.
Since
 \begin{align*}N_EA_{\log T}^- j &= \{(z,r)\in \bH^3:
  T^{-1} <  r\le 1, z\in E\}	
\end{align*}
and $K$ is the stabilizer of $j$ in $G$,
it follows that $$\# [e]\G\cap B_T(E)=
\#\{\gamma (\hat C_0): \gamma (\hat C_0)\text{ contains $(z,r)$ with } z\in E,\;\; T^{-1}< r\le 1\}.   
$$

By the admissibility assumption on $E$,
 we observe that
 $\gamma(\hat C_0)$ contains $(z,r)$ with $z\in E$ and $T^{-1}< r\le 1$
 if and only if the center of $\gamma (C_0)$ lies in $E$ and
the radius of $\gamma( C_0)$ is greater than $T^{-1}$, possibly except for finitely many number (say, $m_0$) of circles.
\end{proof}

\section{Uniform distribution along the family $B_T(E)$ and the Burger-Roblin measure}\label{sbr}
We keep the notations $C_0, H, K, M, A^+, X_0, G,
\{m_x:x\in \bH^3\}$, etc., from section \ref{sectionh}.
Denoting by $dm$ the probability invariant measure on $M$,
\begin{equation} \label{harh} dh=d\mu^{\Leb}_{C_0^\dag}(s)dm\end{equation} is a Haar measure on $H\cong C_0^\dag \times
M$ as $\mu^{\Leb}_{C_0^\dag}$ is $H$-invariant, and
 the following defines a Haar measure on $G$: for $g=ha_rk\in HA^+K$,
\begin{equation}\label{harg} dg=4\sinh r \cdot \cosh r  \; d h dr dm_{j}(k)\end{equation}
where $dm_{j}(k):=dm_{j}(k.X_0^+)$.

We denote by $d\lambda$ the unique measure on $H\ba G$ which is compatible with the choice
of $dg$ and $dh$: for $\psi\in C_c(G)$,
$$\int_G\psi \; dg=\int_{[g]\in H\ba G}\int_{h\in H}\psi(h[g])\;dhd\lambda[g] .$$

For a bounded set $E\subset \c$,
recall that the set $B_T(E)$ in $H\ba G$ is the image of
the set $$KA^+_{\log T}N_{-E}=\{ ka_t n_{-z}\in G:  k\in K, 0\le t<\log T, z\in E\}$$
under the canonical projection $G\to H\ba G$.

The goal of this section is to deduce the following theorem \ref{mtt} from Theorem \ref{os}:
\begin{Thm}\label{mtt} 
Let $\G$ be a non-elementary discrete subgroup of $G$. Suppose that
$|m_\G^{\BMS}| <\infty$ and $\op{sk}_\G(C_0)<\infty$.
Suppose that the natural projection map $\G\cap H\ba \hat C_0 \to \G\ba \bH^3$ is proper.
Then for any bounded Borel subset $E\subset \c$ and for any $\psi\in C_c(\G\ba G)$, we have
$$\lim_{T\to \infty}\frac{1}{T^{\delta_\G}}\int_{g\in B_T(E)}\int_{h\in \G\cap H\ba H}\psi(hg)dhd\lambda(g)
= \frac{\op{sk}_\G(C_0) }{\delta_\G \cdot |m^{\op{BMS}}|}\cdot
 \int_{n\in N_{-E}} m^{\BR}_\G(\psi_n )\; dn  $$
 where $\psi_n\in C_c(\G\ba G)^M$ is given by
$\psi_n(g)=\int_{m\in M}\psi(gmn)dm$ and $dn$ is the Lebesgue measure on $N$.
\end{Thm}

In order to prove this result using Theorem~\ref{os}, it is crucial to understand the shape of the set $B_T(E)$
 in the $HA^+K$ decomposition of $G$. This is one of the important technical steps in the proof. 

\vs
\noindent{\bf{On the shape of $B_T(E)$}:}
Fix a left-invariant metric on $G$. For $\e>0$, let $U_\e$ be the $\e$-ball around $e$ in $G$.
For a subset $W$ of $G$, we set $W_\e=W\cap U_\e$.

\begin{Prop}\label{ksmall}
\begin{enumerate}\item If $a_t\in HKa_sK$ for $s>0$, then $|t|\le s.$
\item Given any $\epsilon>0$, there exists $T_0=T_0(\e)$ such that
$$\{k\in K: a_t k\in HKA^+\text{ for some $t>T_0$}\}\subset K_{\e}M.$$
\end{enumerate}
\end{Prop}
\begin{proof}
Suppose $a_t=hk_1a_sk_2$ for $h\in H, k_1,k_2\in K$.
We note that, as $Aj$ is orthogonal to $\hat C_0$ and $j\in \hat C_0$,
\begin{align*}| t|&= d(\hat C_0, a_t j)=d(\hat C_0, hk_1a_s j)\\
&=d(\hat C_0, k_1a_sj)\le d(j, k_1a_s j)=d(j, a_sj)=s,\end{align*}
proving the first claim.
For the second claim,
suppose $a_tk\in HKa_s$ for some $s\ge 0$. Then
 $ka_{-s}\in a_{-t}HK$.
 Applying both sides to $j\in \bH^3$, 
$k(e^{-s}j)\in a_{-t}\hat C_0$.
 Now $a_{-t}\hat C_0=e^{-t} \hat C_0$ is the northern
hemisphere of Euclidean radius $e^{-t}$ about $0$ in $\bH^3$.

On the other hand $A^-j=(0,1]j$ for $A^-=\{a_{-s}: s\ge 0\}$ and
$K_{\e}\{(0,1]j\}$
 consists of geodesic rays in $\bH^3$ joining $j$ and points
 of $K_{\e} ( 0) \subset \c$. Now $K_{\e}( 0)$ 
contains a disk of radius, say $r_{\e}>0$, centered at $0$ in $\c$,
and hence $K_{\e}\{(0,1]j\}$ contains a Euclidean half ball of radius
 $r_{\e}>0$ centered at $0$ in $\bH^3$. 

Therefore for $t>T_0(\e):=-\log(r_\epsilon)$,
$k(e^{-s}j) \in a_{-t}\hat C_0$ implies that $k(e^{-s}j) \in K_{\e}\{(0,1]j\}$,
in other words, $ka_{-s} K \subset K_{\e} A^-K$. By the uniqueness of
 the left $K$-component, modulo the right multiplication by $M$,
 in the decomposition $G=KA^-K$, it follows that $k\in K_{\e}M$, proving the second claim. 
\end{proof}

%
%
%
%
%

For $t\in \br$ and $T>1$,
set $$K_T(t):=\{k\in K: a_tk\in HKA^+_{\log T} \}.$$
As a consequence of Proposition~\ref{ksmall}, we have the following.

\begin{Cor}\label{property}
\begin{enumerate}
\item For all $0\le t< \log T$, $e\in K_T(t)$.
\item For all $t>\log T$,  $K_T(t)=\emptyset$.
\item For any $\e>0$, there exists $T_0(\e) \ge 1$ such that we have
$$K_T(t)\subset K_\e M\quad\text{for all $t>T_0(\e)$}. $$
\end{enumerate}
\end{Cor}

Thus for any $T>1$,
 \begin{equation}\label{hhkk} 
HKA^+_{\log T}=\cup_{ 0\le t<\log T}Ha_tK_T(t).
\end{equation}

Since $B_T(E)=H\ba HKA^+_{\log T}N_{-E}$,
\eqref{hhkk} together with Corollary \ref{property} shows
that $B_T(E)$ is essentially of the form $H\ba Ha_{\log T}K_\e M N_{-E}$.
The following proposition shows that 
$B_T(E)$ can be basically controlled by the set $H\ba H a_{\log T} N_{-E}$.

\begin{Prop}\label{uk} Fix a bounded subset $E$ of $\c$.
There exists $\ell=\ell(E) \ge 1$ such that
 for all sufficiently small $\e>0$,
 $$a_t km n_z \in H_{\ell \e} m a_t n_z U_{\ell \e} $$
holds for any $m\in M$, $t>0$, $z\in E$,
and  $k\in K_\e$.
\end{Prop}
\begin{proof} Recalling that $N^-$ denotes the lower triangular subgroup of $G$, we note that
the product map $N^-\times A\times M\times N \to G$ is a diffeomorphism at a neighborhood
of $e$, in particular, bi-Lipschitz. Hence
there exists $\ell_1>1$ such that  for all small $\e>0$,
 \begin{equation}\label{k}
 K_{ \e}\subset N_{\ell_1 \e}^{-} A_{\ell_1 \e} M_{\ell_1 \e} N_{\ell_1\e} .\end{equation}

 Similarly due to the $H\times A\times N$ product decomposition of $G_\epsilon$, there exists $\ell_2>1$ such that
 \begin{equation}\label{ue}U_\e\subset H_{\ell_2 \e}A_{\ell_2\e}N_{\ell_2 \e}\end{equation}
for all small $\e>0$ (\cite[Lem 2.4]{GorodnikShahOhIsrael}).
We also have $\ell_3>1$ such that
for all small $\e>0$,
\begin{equation}\label{amn}
A_{(\ell_1+\ell_2) \e}N_{(\ell_1+ \ell_2) \e} M_{\ell_1\e}
\subset U_{\ell _3\e}.
\end{equation}

Now let $t>0, k\in K_\e, m\in M, n\in N$.
Then by \eqref{k}, we may write
 $$k=n^-_1 b_1 m_1 n_1  \in  N_{\ell_1 \e}^{-} A_{\ell_1 \e} M_{\ell_1 \e} N_{\ell_1\e}.$$
Since $a_t n^{-}_1a_{-t}\in N_\e^{-}$ for $t>0$, we
have, by \eqref{ue},
 $$a_t n^{-}_1a_{-t}= h_2b_2m_2n_2\in H_{\ell_2 \e}A_{\ell_2\e}M_{\ell_2 \e}N_{\ell_2 \e} .$$

Therefore
\begin{align*}
a_tk mn&=(a_t n^-_1 a_{-t})(a_t b_1m_1n_1) m n\\&=
(h_2b_2m_2 n_2) a_tb_1m_1n_1m n\\
&=h_2b_2m_2 (a_tb_1 b_1^{-1}a_{-t}) n_2 a_tb_1m_1n_1m n
\\
&= h_2a_t (b_2m_2)b_1(b_1^{-1} a_{-t} n_2 a_t b_1) m_1n_1 mn\\
&\in h_2 a_t A_{(\ell_1 +\ell_2)\e}M_{\ell_2 \e}N_{ (\ell_1+\ell_2) \e} M_{\ell_1\e} mn\quad\text{by \eqref{amn}}
\\ & \subset h_2 a_t U_{\ell_3 \e} m n .\end{align*}

As $E$ is bounded,  there exists $\ell=\ell(E) >\ell_2 $ such that
for all small $\e>0$ and for all $z\in E$,
$$U_{\ell_3 \e} mn_z\subset m n_z U_{\ell\e }  .$$
Since $a_t$ commutes with $m$, we obatin for all $z\in E$ that
$$a_t kmn_z\subset H_{\ell\e} m a_t n_z U_{\ell \e} .$$
\end{proof}

\subsection*{Proof of Theorem \ref{mtt}}
Let $\ell=\ell(E)\ge 1$ be as in Proposition \ref{uk}.
For $\psi\in C_c(\G\ba G)$ and $\e>0$,
we define $\psi_\e^{\pm}\in C_c(\G\ba G)$,
$$\psi_\e^+(g):=\sup_{u\in U_{\ell \e} } \psi(gu)\quad \text{and }\quad \psi_\e^-(g):=\inf_{u\in U_{\ell \e} }\psi(gu) .$$

For a given $\eta>0$,
there exists $\e=\e(\eta)>0$ such that for all $g\in \G\ba G$,
$$|\psi_\e^+(g)-\psi_\e^-(g)|\le \eta$$
 by the uniform continuity of $\psi$.

On the other hand, by Theorem \ref{os},
we have $T_1(\eta)\gg 1$ such that
for all $t>T_1(\eta)$,
\begin{align}\label{pstar}&\int_{h\in \G\cap H\ba H}\psi_\e^+(ha_tn) dh \\ \notag &=
\int_{s\in \G\cap H\ba  C_0^\dagger}\int_{m\in M}\psi_\e^+(s a_t mn)dm d\mu^{\Leb}_{C_0^\dag} (s)\\ \notag &=
(1+O(\eta)) \frac{\op{sk}_\G(C_0)}{|m^{\op{BMS}}|} m^{\op{BR}}_\G(\psi_{\e,n}^+) e^{(\delta_\G -2) t}\end{align}
 where $\psi_{\e,n}^+(g)=\int_{m\in M} \psi_\e^+(gmn) dm$.

As $N_{-E}$ is relatively compact, the implied constant can be taken uniformly over all $n\in N_{-E}$.
Let $T_0(\e)>T_1(\eta)$ be as in
 Proposition \ref{ksmall}.
For $[e]=H\in H\ba G$ and $s>0$,
set $$V_T(s):=\cup_{s\le t < \log T}[e]a_tK_T(t)N_{-E}$$
so that $$B_T(E)=V_T(s)\cup (B_T(E)- V_T(s)).$$

Setting $$\psi^H(g):= \int_{h\in \G\cap H\ba H}\psi(hg) dh ,$$
note that $\psi^H$ is left $H$-invariant as $dh$ is a Haar measure.
We will show that $$\limsup_{T\to\infty}
\frac{1}{T^\delta} \int_{[g]\in V_T(T_0(\e))} \psi^H(g) d\lambda(g) =
(1+O(\eta)) \frac{\op{sk}_\G({C_0})}{\delta_\G \cdot |m^{\op{BMS}}|}\cdot
 \int_{n\in N_{-E}} m^{\BR}_\G(\psi_n ) dn .$$

By Corollary \ref{property}, we have
$$V_T(T_0(\e))\subset \cup_{T_0(\e)\le t\le \log T}[e]a_tK_\e MN_{-E} .$$
 Let $[g]\in V_T(T_0(\e))$, so
 $[g]=[e]a_tkmn$ with $T_0(\e) \le t\le \log T$, $k\in K_\e$, $m\in M$ and $n\in N_{-E}$.
By Proposition \ref{uk}, 
there exist $h_0\in H$ and $ u\in U_{\ell \e }$ such that 
$$a_tkmn =h_0m a_t n u$$
so that $[g]=[e]a_t nu$, since $M\subset H$.

We have
$$\psi^H(g)=\int_{h\in \G\cap H\ba H}\psi(h a_t n u) dh
 \le
\int_{h\in \G\cap H\ba H}\psi_\e^+(ha_t n) dh. $$

The measure $e^{2t}dtdn$ is a right invariant measure of $AN$ and
$[e]AN$ is an open subset in $H\ba G$. Hence $d\lambda(a_tn)$ (restricted to $[e]AN$)
 and $e^{2t}dtdn$
are constant multiples of each other. It follows from the formula of $dg$ that
$d\lambda(a_tn)=e^{2t}dtdn$.
Therefore
$$\int_{[g]\in V_T(T_0(\e))}\psi^H(g) d\lambda(g) \le
\int_{n\in N_{-E}} \int_{T_0(\e)<t \le\log T} \int_{h\in \G\cap H\ba H}
\psi_\e^+(ha_tn)
 dh e^{2t} dtdn .$$

By the choice of $\e=\e(\eta)$, we also have
  $$m^{\BR}_\G(\psi_{\e, n}^+)=(1+O(\eta))m^{\BR}_\G(\psi_n)$$
  where the implied constant depends only on $\psi$.
  Hence by \eqref{pstar},
\begin{multline*}
\int_{n\in N_{- E}}\int_{T_0(\e) < t<\log T}\int_{h\in \G\cap H\ba H}\psi_{\e,n}^+(ha_t) dh e^{2t} dtdn
\\ = (1+O(\eta)) \frac{\op{sk}_{\G}(C_0)}{\delta_\G \cdot |m^{\op{BMS}}|}\cdot
 \int_{n\in N_{-E}} m^{\BR}_\G(\psi_n ) dn  \cdot
  ( T^{\delta_\G} -e^{\delta_\G T_0(\e)}).
\end{multline*}

Hence
\begin{multline*}\limsup_T
\frac{1}{T^{\delta_\G}} \int_{[g]\in V_T(T_0(\e))} \psi^H(g) d\lambda(g)
=
(1+O(\eta)) \frac{\op{sk}_{\G}(C_0)}{\delta_\G \cdot |m^{\op{BMS}}|}\cdot
 \int_{n\in N_{-E}} m^{\BR}_\G(\psi_n ) dn .\end{multline*}

On the other hand, since $\G\ba \G H$ is a closed subset of $\G\ba G$,
 so is $\cup_{0\le t\le s} \G\ba \G H a_t K N_{-\overline E} $ for any fixed $s>0$; in particular,
 its intersection with a compact subset of $\G\ba G$ is compact.

Since
 $$\cup_{[g]\in B_T(E)-V_T(s)}\G\ba \G H g\subset \cup_{0\le t\le s} \G\ba \G H a_t K N_{-\overline E} ,$$
and $\psi$ has compact support, we have, as $T\to \infty$,
$$ \int_{[g]\in B_T(E)-
V_T(T_0(\e))} \int_{h\in \G\cap H\ba H}\psi(hg) dh d\lambda(g) =O(1) .$$
Therefore
$$\limsup_T
\frac{1}{T^\delta}\int_{[g]\in B_T(E)} \psi^H(g) d\lambda(g)
\le (1+O(\eta)) \frac{\op{sk}_{\G}(C_0)}{\delta_\G \cdot |m^{\op{BMS}}|}\cdot
 \int_{n\in N_{-E}} m^{\BR}_\G(\psi_n ) dn .$$
As $\eta>0$ is arbitrary and $\e(\eta)\to 0$ as $\eta\to 0$,
we have
$$\limsup_T
\frac{1}{T^\delta}\int_{[g]\in B_T(E)} \psi^H(g) d\lambda(g)
\le \frac{\op{sk}_{\G}(C_0)}{\delta_\G \cdot |m^{\op{BMS}}|}\cdot
 \int_{n\in N_{-E}} m^{\BR}_\G(\psi_n ) dn .$$
Similarly we can show that
$$\liminf_T
\frac{1}{T^\delta}\int_{[g]\in B_T(E)} \psi^H(g) d\lambda(g)
\ge \frac{\op{sk}_{\G}(C_0)}{\delta_\G \cdot |m^{\op{BMS}}|}\cdot
 \int_{n\in N_{-E}} m^{\BR}_\G(\psi_n ) dn .$$
\qed

\section{On the measure $\omega_{\G}$}\label{smeasure}
In this section we will describe a measure $\omega_{\G}$ on $\c$ and show that the term 
\[
 \int_{n\in N_{-E}} m^{\BR}_\G(\Psi_n )\; dn, 
\]
which appears in the asymptotic expression in Theorem~\ref{mtt}, converges to $\omega_{\G}(E)$ as the support of $\Psi$ shrinks to $[e]$ with $\int_{\G\ba G}\Psi \,dg= 1$.

We keep the notations $ G, K, M, A^+, N, N^- , a_t, n_z, n^-_z$, etc., from section \ref{sectionh}.
Throughout this section,
we assume that $\G$ is a non-elementary discrete subgroup of $G$.
Recall that $\{\nu_x=\nu_{\G, x}:x\in\bH^3\}$ denotes a $\G$-invariant conformal density for $\G$ of dimension $\delta_\G>0$.
\begin{Def}\label{ome} Define a Borel measure $\omega_{\G}$
 on $\c$ as follows: for $\psi\in C_c(\c)$,
 $$\omega_{\G}(\psi)=\int_{z\in \c}e^{\delta_\G\beta_z(x,z+j)}\psi(z) \underline{}d\nu_{\G, x}(z) $$
for $x\in \bH^3$ and $z+j:=(z,1)\in \bH^3$. \end{Def}

In order to see that the definition of $\omega_\G$ is independent of the choice of $x\in \bH^3$, we observe that for any $x_1,x_2\in \bH^3$ and $z\in \c$,
  $$ e^{\delta_\G (\beta_z(x_1, z+j)-\beta_z(x_2, z+j ))}
 \frac{d\nu_{x_1}}{d{\nu}_{x_2}}(z)=
e^{\delta_\G \cdot \beta_z(x_1, x_2)}
 \frac{d\nu_{x_1}}{d{\nu}_{x_2}}(z)= 1 $$
 by the conformality of $\{\nu_{x}\}$.

\begin{Lem}\label{omd} For any $x=p+rj\in \bH^3$ and $\psi\in C_c(\c)$,
$$\omega_\G(\psi):=\int_{z\in \c} (r^{-1}|z-p|^2+r)^{\delta_\G} \psi( z) d\nu_{x}(z) .$$
\end{Lem}
\begin{proof}
It suffices to show that
$$\beta_z( p+rj, z+j)=\log\frac{|z-p|^2+r^2}{r} .$$
We use the fact that the hyperbolic distance $d$ on the upper half space model
of $\bH^3$
satisfies
$$\cosh (d(z_1+r_1j, z_2+r_2j))=\frac{|z_1-z_2|^2+r_1^2+r_2^2}{2r_1r_2}$$
for $z_i+r_i j \in \bH^3$  (cf. \cite{ElstrodtGrunewaldMennickebook}).

Note that \begin{align*}
\beta_z(z, z+j)&=\beta_0(j, -z+p+rj)\\&=
\lim_{t\to \infty} t -d(-z+p+rj, e^{-t}j)
\\ &= \lim_{t\to \infty} t- d(p+rj, z+e^{-t}j) .\end{align*}

Now $$\cosh  d(p+rj, z+e^{-t}j)=\frac{e^t(|z-p|^2+r^2)+e^{-t}}{2r}$$
and hence
$$  e^{d(p+rj, z+e^{-t}j)} +e^{-d(p+rj, z+e^{-t}j)}
= \frac{e^t(|z-p|^2+r^2)+e^{-t}}{r}. $$
Therefore as $t\to \infty$,
$$d(p+rj, z+e^{-t}j)\sim t+\log \frac{|z-p|^2+r^2}{r} .$$
Hence $$\beta_z(p+rj, z+j)=\log\frac{|z-p|^2+r^2}{r} .$$
\end{proof}

\begin{Def}{\rm  For a function $\psi$ on $\c$ with compact support, define
a function $\mathfrak R_\psi$ on $MAN^-N\subset G$ by
$$\mathfrak R_\psi(m a_t n^-_x n_z)=e^{-\delta_\G t} \psi (-z) $$
for $m\in M, t\in \br, x,z\in \c$.
If $\psi$ is the characteristic function of $E\subset \c$,
 we put $\mathfrak R_E=\mathfrak R_{\chi_E}$.
}\end{Def}

Since the product map $M\times A\times N^-\times N\to G$ has a diffeomorphic image,
the above function is well-defined.
\begin{Prop}\label{omd2} For any $\psi\in C_c(\c)$,
 $$\omega_\G(\psi)=\int_{k\in K/M}{\mathfrak R}_\psi (k^{-1}) d\nu_{j}(k(0)).$$
\end{Prop}
\begin{proof} If $k\in K$ with $k^{-1}=ma_tn^-_xn_z \in MAN^-N$,
since $MAN$ fixes $0$,
$$k(0)= n_{-z}(0)=-z .$$

 We note that $\lim_{s\to \infty} a_{-s}(j)=0$
and compute
\begin{align*}
&0=\beta_{-z}(k(j), j)\\& =\beta_{-z}(n_{-z} n^-_{-x} a_{-t} j, j) \\
&=\beta_0( n^-_{-x}a_{-t}j, n_z(j) )\\
& =\lim_{s\to \infty} d( n^-_{-x}a_{-t} j, a_{-s}j)  -d(n_z(j), a_{-s} j)\\
 &=\lim_{s\to \infty} d((a_{s} n^-_{-x}a_{-s}) a_{s-t} j, j) - d(n_z(j), a_{-s} j) \\
  &=\lim_{s\to \infty} d(a_{s-t}j, j)   -d(n_z(j), a_{-s} j)\\
  &= \lim_{s\to \infty} s-t   -d(n_z(j), a_{-s} j)
\end{align*}
and hence
$$-t=\lim_{s\to \infty} d(n_z(j), a_{-s} j) -s=\beta_0( n_z(j), j ) =\beta_{-z}(j, -z+j) .$$

Hence for $k^{-1}\in K\cap MAN^-N$,
$$ {\mathfrak R}_\psi(k^{-1})=e^{\delta_\G \beta_{k(0)}(j, n_{k(0)}(j) )}\psi (k(0)).$$

Since the complement of $NN^-AM/M$ in $K/M$
is a single point and $\nu_j$ is atom-free,
we have
\begin{align*}&\int_{k\in K/M}{\mathfrak R}_E(k^{-1}) d\nu_{j}(k(0))\\
&=\int_{k\in (K\cap NN^-AM)/M}{\mathfrak R}_E(k^{-1}) d\nu_{j}(k(0)) \\
&=\int_{z\in \c}e^{\delta_\G \beta_{k(0)}( j, k(0)+j)} \psi( k(0)) d\nu_{j}(k(0))
\\
&=\int_{z\in \c}
e^{\delta_\G \beta_{-z} (j,-z+j)}\psi(-z) d\nu_{j}(-z)
\\
&=\int_{z\in \c}
e^{\delta_\G \beta_{z} (j, z+j)}\psi(z) d\nu_{j}(z)=\omega_\G(\psi).
\end{align*}
\end{proof}

\begin{Lem} \label{man} If
$(ma_tn_x^-n_z)(m_1 a_{t_1}n_{x_1}^-n_{z_1})=m_0 a_{t_0}n_{x_0}^-n_{z_0}$ in the
$MAN^-N$ coordinates,
then $$t_0=t+t_1+2\log(|1+e^{-t_1}x_1z'|)$$ for some $z'\in \c$ with $|z|=|z'|$.
\end{Lem}
\begin{proof} Note
that if $m_1=\text{diag}(e^{i\theta_1}, e^{-i\theta_1})$, then
$$a_tn_x^-n_zm_1= m_1 a_t n_{e^{i\theta_1}x}^- n_{e^{i\theta_1} z} .$$
Hence we may assume $m_1=m=e$ without loss of generality.
We use the following simple identity for $z,x\in \c$:
\begin{equation}\label{id} n_zn_x^-=\begin{pmatrix}1+xz &0\\ 0&(1+xz)^{-1}\end{pmatrix} n^-_{x(1+xz)}n_{z(1+xz)^{-1}} .
\end{equation}

Hence we have
\begin{align*}&(a_tn_x^-n_z)( a_{t_1}n_{x_1}^-n_{z_1})
\\ &=
(a_{t+t_1})(a_{t_1}^{-1}n_x^-a_{t_1})(a_{t_1}^{-1}n_za_{t_1}) n_{x_1}^-n_{z_1}\\
&=a_{t+t_1} n_{e^{t_1}x}^- n_{e^{-t_1}z} n_{x_1}^-n_{z_1}\\ &=
a_{t+t_1} n_{e^{t_1}x}^-  \begin{pmatrix}1+ e^{-t_1}x_1z &0\\ 0&(1+e^{-t_1}x_1 z)^{-1}\end{pmatrix} n^-_{x_1(1+e^{-t_1}x_1z)} n_{e^{-t_1}z(1+e^{-t_1}zx_1)^{-1}}n_{z_1}
\\ &=m a_{t+t_1+2\log(|1+e^{-t_1}x_1z|)}
  n^-_{x_2 } n_{z_2}
\end{align*}
for appropriate $m\in M$ and $x_2, z_2\in \c$.
\end{proof}

Let $E\subset \c$ be a bounded subset and $U_\e\subset G$ a symmetric $\e$-neighborhood of $e$ in $G$.
For $\e>0$, set
\begin{equation}\label{e}
E_\e^{+}:=U_\e (E)\quad\text{and}\quad E_\e^{-}:=\cap_{u\in U_\e} u(E).\end{equation}

\begin{Lem}\label{casek}\label{ell}  There exists $\ell>0$ such that for all small $\e>0$ and any $g\in U_{\ell \e}$,
$$\int_{k\in K/M}  {\mathfrak R}_E(k^{-1}g) d\nu_{j}(k(0)) =
(1+O(\e))\cdot \omega_\G({E_{\e}^\pm})$$
where the implied constant depends only on $E$.
\end{Lem}
\begin{proof}
Write $k^{-1}=ma_tn_x^-n_z$ and $g=m_1a_{t_1}n_{x_1}^-n_{z_1}\in U_{\e}$.
By Lemma \ref{man}, we have
$k^{-1}g =m_0a_{t_0}n_{x_0}n_{z_0}$ where
$t_0=
t+t_1+2\log(|1+e^{-t_1}x_1z|)$.
Since
 $ {\mathfrak R}_E(k^{-1}g)=e^{-\delta_\G t_0}\chi_{E}(g^{-1}k(0)),
$
we have
\begin{align*}
&\int_{k\in K/M}  {\mathfrak R}_E (k^{-1}g) d\nu_{j}(k(0))\\
&=\int_{k(0)\in g(E)} e^{-\delta_\G t_0} d\nu_j(k(0)) \\
&=\int_{k(0)\in g(E)} e^{-\delta_\G t} e^{-\delta(t_1+2\log (|1+e^{t_1}x_1 z|))} d\nu_j(k(0))\\
&=(1+O(\e)) \int_{k(0)\in E_{\e}^\pm} e^{-\delta_\G t} e^{-\delta(t_1+2\log (|1+e^{t_1}x_1 z|))} d\nu_j(k(0)).
\end{align*}

Since $t_1=O(\e), x_1=O(\e)$ and $z=-k(0)\in -g(E)\subset -E_\e^+$,
$$t_1+2\log (|1+e^{-t_1}x_1 z|)=O(\e)$$
where the implied constant depends only on $E$.
Hence
\begin{align*}&
\int_{k\in K/M}  {\mathfrak R}_E (k^{-1}g) d\nu_{j}(k(0))
\\&=(1+O(\e)) \int_{k(0)\in E_{\e}^\pm} e^{-\delta_\G t} d\nu_j(k(0))\\
&=(1+O(\e)) \int_{k\in K}  {\mathfrak R}_{E_{\e}^\pm}(k^{-1}) d\nu_{j}(k(0))\\ &=(1+O(\e))\cdot \omega_\G
({E_{\e}^\pm}) .\end{align*}
\end{proof}

For $\e>0$,
let $\psi^\e$ be a non-negative continuous function in $C(G)$ with support in $U_{ \e}$
with integral one and
 $\Psi^\e\in C_c(\G\ba G)$ be the $\G$-average of $\psi^\e$:
 $$\Psi^\e(\G g):=\sum_{\gamma\in \G}\psi^\e(\gamma g) .$$

We define $\Psi_E^{\e} \in C_c(\G\ba G)^M$ by
$$\Psi_{E}^{\e}(g):=\int_{z\in -E}\int_{m\in M}\Psi^\e(gmn_{z}) dmdz .$$

\begin{Lem}\label{ccc} For a bounded Borel subset $E\subset \c$,
 there exists $c=c(E)>1$ such that for all small $\e>0$,
$$ (1-c\cdot \e)\cdot
\omega_{\G}( E_\e^-)\le
  m^{\BR}_\G(\Psi^\e_{E} )  \le (1+c\cdot \e)\cdot
\omega_{\G}( E_\e^+)  . $$

\end{Lem}

\begin{proof} Note that 
$N^-$ is the expanding horospherical subgroup for the right action of $a_t$, i.e.,
$N^-=\{g\in G: a_t g a_{-t}\to e\text{ as $t\to \infty$}\}$.
We have for $\psi\in C_c(G)^M$,
$$\tilde m^{\BR}_\G(\psi)=\int_{KAN^-}\psi(ka_t n^-) e^{-\delta_\G t} dn dt d\nu_{j}(k(0)) $$
(cf. \cite[6.2]{OhShahGFH}).
We note that $d(a_tn_x^-m n_z)=dtdxdm dz$ is the restriction of the
Haar measure $dg$ to $AN^- N\subset G/M$.

We deduce
\begin{align*}& m^{\BR}_\G(\Psi_E^\e)=
\int_{z\in {- E}} \tilde m^{\BR}(\psi^\e_{n_z} ) dz\\
 &=\int_{ z\in {-E}} \int_{KAN^-} \int_{m\in M}
\psi^\e(ka_tn_x^- m n_z) e^{ - \delta_\G t}dm dxdtd\nu_{j}(k(0)) dz \\
&=\int_{k\in K}\int_{AN^-MN}\psi^\e(k(a_tn^-_x m n_z)) \chi_{-E}(z) e^{ - \delta_\G t} dxdt dm dz d\nu_{j}(k(0)) \\
&=\int_{k\in K}\int_{g\in G } \psi^\e(k g) {\mathfrak R}_E(g) dg d\nu_{j}(k(0))\\
&=\int_{g\in U_\e} \psi^\e(g) \left(\int_{k\in K}  {\mathfrak R}_E(k^{-1}g) d\nu_{j}(k(0)) \right) dg. \end{align*}
Hence by Lemma \ref{ell} and the identity $\int_{U_\e} \psi^\e dg=1$,
 we have
$$m^{\BR}_\G(\Psi_E^\e)=
  (1+O(\e))\omega_{\G}(E_{\e}^{\pm}) .$$
\end{proof}

\begin{Cor}\label{inter}
If $\omega_\G(\partial(E))=0$, then
$$\omega_\G(E)=\lim_{\e\to 0} m^{\BR}_\G(\Psi^\e_{E} ) .$$
\end{Cor}
\begin{proof}
For any $\eta>0$, there exists $\e=\e(\eta)$ such that
$\omega_{\G}(E_\e^+-E_{\e}^-)<\eta$.

Together with  Lemma \ref{ccc}, it implies that
$$  m^{\BR}_\G(\Psi^\e_{E} ) = (1+O(\e))(1+O(\eta))\omega_\G(E)$$
and hence the claim follows.
\end{proof}

\section{Conclusion: Counting circles}\label{conclusion}
Let $\G<G:=\PSL_2(\c)$ be a non-elementary discrete group
with $|m^{\BMS}_\G|<\infty$.
 Suppose that $\mathcal P:=\G(C)$ is a locally finite circle packing.

Recall that $$\op{sk}_{\G}(\P)=\op{sk}_{\G}(C):=
 \int_{s\in \op{Stab}_{\G} (C^\dagger)\ba C^\dagger}  e^{\delta_\G
\beta_{s^+}(x,s)}d\nu_{\G, x}(s^+), $$
 where $C^\dagger$ is the set of unit normal vectors to $\hat C$.
It follows from the conformal property of $\{\nu_{\G, x}\}$
that $\op{sk}_{\G}(C)$
is independent of the choice of $C\in\Gamma(C)$, and hence
is an invariant of the packing $\G(C)$.

Theorem \ref{m1} is an immediate consequence of the following statement.

\begin{Thm}\label{mmt2} Suppose that $\op{sk}_\G(C)<\infty$.
For any bounded Borel subset $E$ of $\c$ with $\omega_\G(\partial(E))=0$,
 we have
\begin{equation}\label{nt}
\lim_{T\to \infty}\frac{N_T(\mathcal P, E)}{T^{\delta_\G} }= 
\frac{\op{sk}_{\G}(\P)}{\delta_\G \cdot |m^{\BMS}_\G|} \cdot \omega_{\G}(E).\end{equation}
Moreover $\op{sk}_\G(C) >0$ if $\P$ is infinite.
\end{Thm}
The second claim on the positivity of $\op{sk}_\G(C)$ follows from the second claim
of Theorem \ref{os} and Lemma \ref{infinite}.

We will first prove Theorem \ref{mmt2} for the case when $C$ is the unit circle $C_0$ centered
at the origin and deduce the general case from that.

\subsection*{The case of $C=C_0$.}
Fix $\eta>0$.
As  $\omega_\G(\partial (E))=0$,
there exists $\e=\e(\eta)>0$ such that
\begin{equation}\label{ej}\omega_\G( {E_{4\e}^+}  - E_{4\e}^-)\le \eta \end{equation}
where $E_{4\e}^{\pm}$ is defined as in \eqref{e}:
$E_{4\e}^{+}:=U_{4\e} (E)$ and $ E_{4\e}^{-}:=\cap_{u\in U_{4\e}} u(E)$.

We can find a $\P$-admissible Borel subset $\tilde E_{\e}^{+}$ such that
$E\subset \tilde E_{\e}^{+}\subset E_{\e}^+$  by adding
all the open disks inside $E_{\e}^{+}$
 intersecting the boundary of $E$.
  Similarly we can find a $\P$-admissible Borel subset $\tilde E_{\e}^{-}$ such that
 $E_{\e}^-\subset\tilde E_{\e}^-\subset E$ by adding all the open disks inside $E$
 intersecting the boundary of $E_{\e}^{-}$.
 By the local finiteness of $\P$, there are only finitely many circles intersecting $E$ (resp.
 $\tilde E_{\e}^{-}$) which are not contained in $\tilde E_{\e}^{+}$ (resp. $E$). Therefore
there exists $q_\e\ge 1$ (independent of $T$) such that
\begin{equation}\label{npt}
N_T(\P, \tilde E_{\e}^{-} ) - q_\e\le N_T(\P, E) \le N_T(\P, \tilde E_{\e}^{+} ) +q_\e .\end{equation}


Recalling the set $B_T(\tilde E^{\pm}_\e)= H\ba H KA^+_{\log T}N_{-\tilde E^{\pm}_\e}\subset H\ba G ,$
it follows from Proposition \ref{tran} and \eqref{npt} that for all $T\gg 1$,
\begin{equation}\label{tft}
\# [e]\G\cap B_T(\tilde E^{-}_\e)-m_0 \le N_T(\G (C_0), E)\le \# [e]\G\cap B_T(\tilde E^{+}_\e) +m_0 \end{equation}
for some fixed $m_0=m_0(\e)\ge 1$.


\begin{Lem}\label{strrr}
 There exists $\ell >0$ such that for all $T> 1$ and for all small $\e>0$,
 $$KA_{\log T}^+U_{\e}\subset K A^+_{\log T +\e}N_{\ell \e}$$
where $N_{\ell \e}$ is the $\ell \e$-neighborhood of $e$ in $N$.
\end{Lem}
\begin{proof}
We may write $U_\e=M_\e N^-_\e A_\e N_\e=K_\e A_\e N_\e$ up to uniform Lipschitz constants.
For $u=mn^-a n\in M_\e N^-_\e A_\e N_\e$,
$a_t u= m (a_tn^-a_{-t}) a_ta n$.
Since $a_t n^-a_{-t}\in U_{\e}$ for $t>0$,
we may write it as $k_1 a_1 n_1\in K_{\e} A_{\e} N_{ \e}$.
Hence for $0<t<\log T$, we have
 $(a^{-1} a_{-t} n_1 a_t a) \in N_\e$ and
$$a_tu =(mk_1) (a_1 a_ta) (a^{-1} a_{-t} n_1 a_t a) 
 n \in K A^+_{\log T + 2\e }N_{2\e} .$$
This proves the claim.
\end{proof}

\begin{Lem}[Stability of $KAN$-decomposition]\label{strong}\label{strong2} There exists $\ell_0>0$
 (depending on $E$) such that for all $T> 1$ and for all
small $\e>0$,
\begin{equation*}
KA_{\log T}^+N_{-\tilde E^+_\e} U_{\ell_0\e} \subset KA_{\log T +\e}^+  N_{-E_{2\e}^+} ;\end{equation*}
\begin{equation*}
KA_{\log T-\e}^+N_{-E_{2\e}^-} \subset   KA_{\log T}^+ (\cap_{u\in U_{\ell_0\e}} N_{-\tilde E^-_\e} u) .\end{equation*}
\end{Lem}

\begin{proof}
There exists $\ell_0>0$ depending on $E$
 such that $N_{-\tilde E^+_\e} U_{\ell_0\e}\subset U_\e N_{-E_{2\e}^+}$.
Hence the first claim follows from Lemma \ref{strrr}. The second claim can be proved similarly.
\end{proof}




For $\e>0$, define functions $F_T^{\e, \pm}$ on $\G\ba G$:
 $$F_T^{\e, +}(g):=\sum_{\gamma\in (H\cap \G)\ba \G}\chi_{B_{e^{\e}T}(N_{-E_{2\e}^+})}([e]\gamma g);\quad
 F_T^{\e, -}(g):=\sum_{\gamma\in (H\cap \G)\ba \G}\chi_{B_{e^{-\e}T}(N_{-E_{2\e}^-})}([e]\gamma g) .$$

Let $\ell_0>0$ be as in Lemma \ref{strong}.
Without loss of generality, we may assume that $\ell_0<\ell$ 
for $\ell$ as in Lemma \ref{ell}.

\begin{Lem} For all $g\in U_{\ell_0 \e }$ and $T\gg 1$,
\begin{equation}\label{ft} F_T^{\e, +}(g)-m_0
\le  N_T(\G (C_0), E) \le  F_T^{\e, +}(g)+m_0.
\end{equation}\end{Lem}
\begin{proof} Note that, since $U_{\ell_0\e}$ is symmetric, for any $g\in  U_{\ell_0 \e }$,
$$\# [e]\G\cap B_T(\tilde E^+_\e) \le \# [e]\G\cap B_T(\tilde E^+_\e) U_{\ell_0 \e} g^{-1} \le
\# [e]\G g \cap B_{e^{\e}T}(N_{-E_{2\e}^+}) ,$$
by Lemma \ref{strong}, which proves the second inequality by \eqref{tft}.
The other
 inequality can be proved similarly.
\end{proof}

For $\e>0$,
let $\psi^\e$ be a non-negative continuous function in $C(G)$ with support in $U_{\ell_0 \e}$
with integral one and
 $\Psi^\e\in C_c(\G\ba G)$ be the $\G$-average of $\psi^\e$:
 $$\Psi^\e(\G g):=\sum_{\gamma\in \G}\psi^\e(\gamma g) .$$

By integrating \eqref{ft} against $\Psi^\e$, we have
$$\la F_{T}^{\e,-}, \Psi^\e\ra -m_0
\le  N_T(\G (C_0), E) \le \la F_{T}^{\e,+}, \Psi^\e\ra +m_0. $$

Since \begin{align*}
\la F_{T}^{\e,+}, \Psi^\e \ra \notag  &=\int_{\G\ba G}
 \sum_{\gamma\in \G\cap H\ba \G}\chi_{B_{e^\e T}(N_{-E_{2\e}^+})}([e]\gamma g )
\Psi^\e(g) \; dg\\ &= \int_{g\in \G\cap H\ba G} \chi_{B_{e^\e T}(N_{-E_{2\e}^+})} ([e] g) \Psi^\e(g) \; dg \notag
\\ &=\int_{[g]\in B_{e^\e T}(N_{-E_{2\e}^+})}\int_{h\in \G\cap H\ba H}\Psi^\e (hg)\, dhd\lambda(g)
\end{align*}
we deduce from Theorem \ref{mtt} and Lemma \ref{locfinite}
that
\begin{equation}\label{sim2} \la F_{T}^{\e,+}, \Psi^\e \ra \sim
\frac{\op{sk}_{\G}(C_0)}{\delta_\G \cdot |m^{\op{BMS}}_\G|}\cdot
 \int_{n\in N_{-E_{2\e}^+} } m^{\BR}_\G(\Psi_n^\e ) dn  \cdot T^{\delta_\G} \cdot e^{\e \delta_\G} \end{equation}
where $\Psi_n^\e (g)=\int_{m\in M}\Psi^\e(gmn)dm$.

 Therefore by applying Lemma \ref{ccc} to \eqref{sim2} and using \eqref{ej}, we deduce
\begin{align*} \limsup_T \frac{\la F_{T}^{\e,+}, \Psi^\e \ra}{T^{\delta_\G}}
 &\le (1+\e) \frac{\op{sk}_{\G}(C_0)}{\delta_\G \cdot |m^{\op{BMS}}_\G|}\cdot
 \int_{n\in N_{-E_{2\e}^+} } m^{\BR}_\G(\Psi_n^\e ) dn  \\ & \le
 (1+\e)(1+c\e) \frac{\op{sk}_{\G}(C_0)}{\delta_\G \cdot |m^{\op{BMS}}_\G|}\cdot
 \omega_{\G}(E_{4\e}^+) \\&
 \le (1+ c_1 \eta )(1+ c_2\e ) \frac{\op{sk}_{\G}(\G(C_0))}{\delta_\G \cdot |m^{\op{BMS}}_\G|}\cdot
 \omega_{\G}(E) \end{align*}
where the constants $c,c_1, c_2$ depend only on $E$.

Similarly, we have
\begin{align*} \liminf_T \frac{\la F_{T}^{\e,+}, \Psi^\e \ra}{T^{\delta_\G}}
\ge (1-c_1\eta)(1-c_2\e) \frac{\op{sk}_{\G}(C_0)}{\delta_\G \cdot |m^{\op{BMS}}_\G|}\cdot
 \omega_{\G}(E).
\end{align*}

As $\eta>0$ is arbitrary and $\e=\e(\eta)\to 0$ as $\eta\to 0$, we have
$$ 
\lim_{T\to \infty}\frac{N_T(\G (C_0), E)}{T^{\delta_\G} }= 
\frac{\op{sk}_{\G}(C_0)}{\delta_\G \cdot |m^{\BMS}_\G|} \cdot \omega_{\G}(E).$$

This proves Theorem \ref{mmt2} for $C=C_0$.

\subsection*{The general case}
 Let $r>0$ be the radius of $C$ and $p\in \c$ the center of $C$.
 Set
  $$g_0=n_{p}a_{\log r}=\begin{pmatrix} 1& p\\ 0 &1\end{pmatrix}\begin{pmatrix} \sqrt {r}& 0\\ 0 &\sqrt{r^{-1}} \end{pmatrix} .$$

 Then $g_0^{-1}(z)= r^{-1}(z-p)$ for $z\in \c$ and $g_0^{-1}(C)=C_0 .$

Setting $\G_0=g^{-1}_0\G g_0$, we have  \begin{align*}
N_T(\Gamma (C), E)&=\#\{C\in \Gamma(g_0(C_0)): C^\circ \cap E\ne\emptyset, \op{Curv}(C)<T\}\\
&=\#\{g_0^{-1}(C)\in \Gamma_0(C_0): C^\circ \cap E\ne\emptyset, \op{Curv}(C)<T\}\\
&=\#\{C_*\in \Gamma_0(C_0): C_*^\circ \cap g_0^{-1} (E)\ne\emptyset, \op{Curv}(C_*)<r^{-1} T\}\\
&=N_{r^{-1}T}(\Gamma_0 (C_0), g_0^{-1}(E)).
\end{align*}

We claim that
\begin{equation}\label{haha}
\frac{1} { |m_{\G_0}^{\BMS}|}  \cdot
\op{sk}_{\G_0}(\G_0(C_0)) \cdot r^{-\delta_\G}\cdot \omega_{\G_0} (g^{-1}_0(E))=
\frac{1} { |m_{\G}^{\BMS}|}  \cdot
\op{sk}_{\G}(\G(C)) \cdot  \omega_{\G} (E).
\end{equation}
Note that
the each side of the above is independent of the choices of conformal densities of $\G_0$ and $\G$ respectively.

Fixing a $\G$-invariant conformal density $\{\nu_{\G, x}\}$
of dimension $\delta_\G$,
set $$\nu_{\G_0, x}:={g_0}^*\nu_{\G, g_0(x)}$$
where $g_0^*\nu_{\G, g_0(x)}(R)=\nu_{\G, g_0(x)}(g_0(R)).$
It is easy to check that $\nu_{\G_0,x}$ is supported on  $\Lambda(\G_0)=g_0\Lambda(\G)$
and satisfies
$$\frac{d\nu_{\G_0,x}}{d\nu_{\G_0,y}}(z)=e^{-\delta_\G\beta_z(x,y)};\quad
 \gamma_*\nu_{\G_0, x}=\nu_{\G_0,\gamma(x)}$$
for all $x,y\in \bH^3$, $\gamma\in \G_0$ and $z\in \hat \c$.

 Hence $\{\nu_{\G_0, x}:x\in \bH^3\}$ is a $\G_0$-invariant conformal density
of dimension $\delta_\G=\delta_{\G_0}$ and satisfies that for $f\in C_c(\c)$
$$\int_{g_0(z)\in E} f(z) d\nu_{ \G_0,x}(z)=
\int_{z\in E} f(g^{-1}_0(z)) d\nu_{\G,g_0(x)}(z) .$$

We consider the Bowen-Margulis-Sullivan measures  $m^{\BMS}_{\G}$ and
$m^{\BMS}_{\G_0}$
on $\G\ba \T^1(\bH^3)$ and $\G_0\ba \T^1(\bH^3)$  associated to $\{\nu_{\G, x}\}$ and
$\{\nu_{\G_0, x}\}$, respectively.

\begin{Lem}\label{fione} For a bounded Borel function $\psi$ on $\G\ba \T^1{(\bH^3)}$,
consider a function $\psi_{g_0}$ on $\G_0\ba \T^1(\bH^3)$  given by $\psi_{g_0}(u):= \psi(g_0 (u)).$
Then
$$m_{ \G_0 }^{\BMS}(\psi_{g_0})=m_{\G}^{\BMS}(\psi) .$$
In particular,
$|m_{\G_0}^{\BMS}|=|m_{\G}^{\BMS}| .$
\end{Lem}
\begin{proof} Note that
if $v=g(u)$, then
$$ \beta_{u^{\pm}}(x, \pi(u))= \beta_{v^{\pm}}(g(x), \pi(v)).$$
Since $\nu_{\G_0, x}=g_0^*\nu_{\G, g_0(x)}$, we have
\begin{align*} &m_{\G_0 }^{\BMS}(\psi_{g_0})\\ &
=\int_{u\in \G_0\ba \T^1(\bH^n)} \psi(g_0(u))e^{\delta_\G \beta_{u^+}(x, \pi(u))}\;
 e^{\delta_\G \beta_{u^-}(x,\pi(u)) }\;
d\nu_{\G_0,x}(u^+) d\nu_{\G_0 ,x}(u^-) dt\\ &=
\int_{v\in \G\ba \T^1(\bH^n)} \psi(v)e^{\delta_\G \beta_{v^+}(g_0(x), \pi(v))}\;
 e^{\delta_\G \beta_{v^-}(g_0(x),\pi(v)) }\;
d\nu_{\G, g_0(x)}(v^+) d\nu_{\G, g_0(x)}(v^-) dt\\ & =
m^{\BMS}_{\G}(\psi) .\end{align*}
\end{proof}

Similarly, we can verify:
\begin{Lem} \label{fitwo} For any $x\in \bH^3$,
 $$\int_{s\in \op{Stab}_{\G_0}(C_0^\dagger)\ba C_0^\dagger}e^{\delta_\G \beta_{s^+}(x, s)} d\nu_{ \G_0, x}(s^+)
= \int_{s\in \op{Stab}_{\G}(C^\dagger)\ba C^\dagger}e^{\delta_\G \beta_{s^+}(g_0(x), s)} d\nu_{\G, g_0(x) }(s^+);$$
that is,
$\op{sk}_{\G} (\G(C))=\op{sk}_{\G_0} (\G_0(C_0)) .$
\end{Lem}

\begin{Lem} For any bounded Borel subset $E\subset \c$,
$$\omega_{\G_0}(g^{-1}_0(E)) = r^{\delta_\G} \omega_{\G}(E) .$$
\end{Lem}
\begin{proof} Since $g_0^{-1}(z)=r^{-1}(z-p)$,
 $r$ is the linear distortion of the map $g_0^{-1}$ in the Euclidean metric, that is,
$r=\lim_{w\to w_0}\frac{|g_0^{-1}(w)-g_0^{-1}(w_0)|}{|w-w_0|}$ for any $w_0\in \c$.
Hence
 $$d\nu_{\G, g_0(j)}(w)=r^{\delta_\G} \frac{(|w|^2+1)^{\delta_\G}}{(|g_0^{-1}(w)|^2+1)^{\delta_\G}}  d\nu_{\G, j}(w).$$

Since  $\nu_{\G_0, x}=g_0^*\nu_{\G, g_0(x)}$, we deduce
\begin{align*}\omega_{\G_0}(g^{-1}_0(E))&=
\int_{z\in g_0^{-1}(E)}( |z|^2+1)^{\delta_\G} d\nu_{\G_0, j}(z)\\&=
\int_{u\in E} (|g_0^{-1}(u)|^2+1)^{\delta_\G} d\nu_{\G, g_0(j)}(u)
\\ &=r^{\delta_\G} \int_{ u\in E }(|u|^2+1)^{\delta_\G} d\nu_{\G,j}(u)
\\ &=r^{\delta_\G} \omega_\G(E).
\end{align*}
\end{proof}

This concludes a proof of \eqref{haha}.
Therefore, since $\sk_{\G_0}(C_0)<\infty$ and $|m_{\G_0}^{\BMS}|<\infty$, the previous case of $C=C_0$ yields that
\begin{align*}
\lim_{T\to \infty}\frac{1}{T^{\delta_\G}} N_T(\Gamma (C), E)
&=\lim_{T\to \infty}\frac{1}{T^{\delta_\G}} 
N_{r^{-1}T}(\Gamma_0 (C_0), g_0^{-1}(E))\\ & =
 \frac{1}{\delta_{\G_0} \cdot |m_{\G_0}^{\BMS}|}  \cdot
\op{sk}_{\G_0}(C_0) \cdot r^{-\delta_\G}\cdot \omega_{\G_0} (g^{-1}_0(E))
\\& = \frac{1}{\delta_\G \cdot |m_{\G}^{\BMS}|}  \cdot
\op{sk}_{\G} (C) \cdot \omega_{\G}(E).\end{align*}
This completes the proof of Theorem \ref{mmt2}. \qed

\bibliographystyle{plain}

\end{document}